\RequirePackage{etoolbox}
\csdef{input@path}{{style/}{graphics/}}
\documentclass[MSNbibl,number,citesort,seceqn,dvips]{arxbj}
\usepackage{upgreek}
\usepackage{mathrsfs}

\aid{0}
\volume{21}
\issue{1}
\pubyear{2015}
\firstpage{303}
\lastpage{334}
\doi{10.3150/13-BEJ568} 

\makeatletter

\newcommand{\rright}{\right}
\newcommand{\lleft}{\left}
\newtheorem{theorem}{Theorem}[section]
\newremark{rem}{Remark}[section]
\newtheorem{lem}{Lemma}[section]
\newproclaim{defn}{Definition}[section]
\newproclaim{assum}[theorem]{Assumption}
\renewcommand{\pi}{\uppi}
\def\a{\alpha}
\def\e{\varepsilon}
\def\l{\lambda}
\def\n{\nu}
\def\si{\sigma}
\def\t{\tau}
\def\f{\varphi}
\def\h{\widehat}
\def\f{\phi}

\def\bE{\mathbf{E}}

\def\cD{\mathcal{D}}
\def\cF{\mathcal{F}}
\def\cU{\mathcal{U}}
\def\cX{\mathcal{X}}

\def\hC{\mathbb{C}}
\def\hD{\mathbb{D}}
\def\hE{\mathbb{E}}
\def\hF{\mathbb{F}}
\def\hN{\mathbb{N}}
\def\hP{\mathbb{P}}
\def\hR{\mathbb{R}}

\def\sB{\mathscr{B}}
\def\sE{\mathscr{E}}
\def\sH{\mathscr{H}}

\def\b1{\mathbf{1}}

\def\sfrac#1#2{#1/#2}

\makeatother

\begin{document}
\begin{frontmatter}

\title{Stochastic differential equations driven by fractional Brownian
motion and Poisson point process}
\runtitle{SDE driven by FBM and Poisson point process}

\begin{aug}
\author[a]{\inits{L.}\fnms{Lihua} \snm{Bai}\thanksref{a}\ead[label=e1]{lhbai@nankai.edu.cn}} 
\and
\author[b]{\inits{J.}\fnms{Jin} \snm{Ma}\corref{}\thanksref{b}\ead[label=e2]{jinma@usc.edu}}
\address[a]{Department of Mathematical Sciences, Nankai University,
Tianjin 300071, China.\\
\printead{e1}}

\address[b]{Department of Mathematics,
University of Southern California, Los Angeles, CA 90089, USA.\\
\printead{e2}}
\end{aug}

\received{\smonth{6} \syear{2012}}
\revised{\smonth{6} \syear{2013}}

%
\begin{abstract}
In this paper, we study a class of stochastic differential equations
with additive noise that contains a fractional Brownian motion (fBM)
and a Poisson point process of class (QL). The differential equation of
this kind is motivated by the reserve processes in a general insurance
model, in which the long term dependence between the claim payment and
the past history of liability becomes the main focus. We establish some
new fractional calculus on the fractional Wiener--Poisson space, from
which we define the weak solution of the SDE and prove its existence
and uniqueness. Using an extended form of Krylov-type estimate for the
combined noise of fBM and compound Poisson, we prove the existence of
the strong solution, along the lines of Gy\"ongy and Pardoux (\textit{Probab. Theory Related Fields}
\textbf{94} (1993) 413--425). Our result
in particular extends the one by Mishura and Nualart (\textit
{Statist. Probab. Lett.} \textbf{70} (2004) 253--261).
\end{abstract}

%
\begin{keyword}
\kwd{discontinuous fractional calculus}
\kwd{fractional Brownian motion}
\kwd{fractional Wiener--Poisson space}
\kwd{Krylov estimates}
\kwd{Poisson point process}
\kwd{stochastic differential equations}
\end{keyword}

\end{frontmatter}

\section{Introduction}
\setcounter{equation}{0}

In this paper, we are interested in the following stochastic
differential equation (SDE):
%
\begin{eqnarray}
\label{SDE0} 
X_{t}=x+\int_{0}^{t}b(s,X_{s})\,\mathrm{d}s+
\si B_{t}^{H}-L_{t}, \qquad   t\in[0,T],
\end{eqnarray}
where $B^H=\{B^H_t\dvtx t\ge0\}$ is a \textit{fractional Brownian motion}
with Hurst parameter $H\in(0,1)$, defined on a given filtered
probability space $(\Omega, \cF, \hP; \hF)$, with $\hF=\{\cF
_t\dvtx t\ge
0\}$ being a filtration that satisfies the \textit{usual hypotheses}
(cf., e.g., \cite{prot}); and $L=\{L_t\dvtx  t\ge0\}$ is a 
Poisson point process of class (QL), independent of $B^{H}$. More
precisely, we assume that $L$ takes the form
%
\begin{eqnarray}
\label{L} L_t=\int_0^t\int
_{\hR}f(s,x)N_p(\mathrm{d}s,\mathrm{d}x), 
 \qquad  t\ge0,
\end{eqnarray}
where $f$ is a deterministic function, and $p$ is a stationary Poisson
point process whose counting measure $N_p$ is a Poisson random measure
with L\'evy measure $\n$ (see Section~\ref{sec2} for more details).

One of the motivations for our study is to consider a general reserve
process of an insurance company, perturbed by an additive noise that
has long term dependency.
A~commonly seen perturbed reserve (or surplus) model is of the
following form:
%
\begin{eqnarray}
\label{U} U_{t}=x+c(1+\rho)t+\e W_{t}-L_{t},\qquad
  t\in[0,T].
\end{eqnarray}
Here $x\geq0$ denotes the initial surplus, $c>0$ is the premium rate,
$\rho>0$ is the
``safety'' (or expense) loading, $\e>0$ is the perturbation parameter,
$W=\{W_{t}\dvtx t\ge0\}$ is a Brownian motion, which represents
an additional uncertainty coming from either the aggregated claims or
the premium income, $L_t$ denotes cumulated claims up to time $t$, and
finally, $T>0$ is a fixed time horizon.
We refer the reader to the well-referred book \cite{RSST}, Chapter~13,
and the references therein for more explanations of such models.

In this paper, we are particularly interested in the case where the
diffusion perturbation
term possesses long-range dependence. Such a phenomenon has been noted
in insurance models based on the observations that the claims often
display long memories due to extreme weather, natural disasters, and
also noted in casualty insurance such as automobile third-party
liability (cf. e.g., \cite{Deb,frango1,frango2,frango3,HP,mich1,mich2} and references therein). A~reasonable refinement that reflects
the long memory but also
retains the original features of the aggregated claims is to assume
that the Brownian motion
$W$ in (\ref{U}) is replaced by a fractional Brownian motion $B^H$,
for a certain Hurst
parameter $H\in(0,1)$. In fact, if we assume further that in addition
to the premium income, the company also receives interest of its
reserves at time with interest rate $r>0$, and that the safety loading
$\rho$ also depends on the current reserve value, one can argue that
the reserve process $X$ should satisfy an SDE of the form of (\ref
{SDE0}) with
\[
b(t,x)=rx+c\bigl(1+\rho(t,x)\bigr), \qquad  (t,x)\in[0,T]\times\hR.
\]


The main purpose of this paper is
to find the minimum conditions on the function $b$ under which the SDE
(\ref{SDE0}) is well posed, in both weak and strong
sense.
In the case when $L\equiv0$, the SDE (\ref{SDE0}) becomes
one driven by an (additive) fBM and the similar issues were
investigated by Nualart and Ouknine
\cite{NO} and Hu, Nualart and Song \cite{HN}. One of the main results is
that, unlike the ordinary differential equation
case, the well-posedness of the SDE can be established under only some
integrability conditions, and in particular, no Lipschitz continuity is
required for uniqueness. The main
idea is to use a Krylov-type estimate to obtain a comparison theorem,
whence the pathwise uniqueness. Such a scheme was utilized by Gy\"
ongy and Pardoux \cite{GyongyPard} when studying the quasi-linear SPDEs,
and has been a frequently
used tool to treat the SDEs with non-Lipschitz coefficients, as an
alternative to the well-known Yamada--Watanabe theorem.
In fact, this method is even more crucial in the current case, as the
usual Yamada--Watanabe theorem type of argument does not seem to work
due to the lack of independent increment property of an fBM.

The main difficulty in the study of SDE (\ref{SDE0}), however, is the
presence of the jumps. In the case when $H>1/2$, Mishura and Nualart
\cite{MN} studied the existence of weak solution of SDE (\ref{SDE0})
with $L\equiv0$, and the coefficient $b$
is allowed to have finitely many discontinuities in its spatial
variable $x$. By a simple transformation (e.g., setting
$\tilde X=X-L$), our result in a sense extends their result to a more
general case in which $b$ possesses countably many discontinuities in
$x$. More importantly, we remove the extra assumption that $H<(1+\sqrt {5})/4$ in \cite{MN} when the number of jumps is finite. To our best
knowledge, the fractional calculus applying to SDE driven by both fBM
and Poisson point process is
new.


The rest of the paper is organized as follows. In Section~\ref{sec2}, we review
briefly the
basics on fBM and some fractional calculus that is needed in this
paper. In Section~\ref{sec3}, we
prove a Girsanov theorem and in Section~\ref{sec4} we apply it to study the
existence of the weak
solution. In Section~\ref{sec5}, we address the uniqueness issue, in both weak
and strong forms,
and in Section~\ref{sec6} we study the existence of the strong solution.


\section{Preliminaries}\label{sec2}

In this section, we review some of the basic concepts in fractional
calculus and
introduce the notion of (canonical) fractional Wiener--Poisson spaces
which will be
the basis of our study. Throughout this paper, we denote $\bE$ (also
$\bE_1, \ldots$) for a generic
Euclidean space, whose inner products and norms will be denoted as
the same ones $\langle\cdot,\cdot\rangle$ and $|\cdot|$,
respectively; and denote $\|\cdot\|$ to be
the norm of a generic Banach space. Let $\cU\subset\bE$ be a
measurable subset. We shall
denote by $L^p(\cU;\bE_1)$, $0\le p<\infty$, the space of all $\bE
_1$-valued measurable function $\f(\cdot)$ defined on $\cU$ such
that $\int_{\cU}|\f(t)|^p\,\mathrm{d}t<\infty$ ($p=0$ means merely
measurable). For each $n\in\mathbb{N}$, $\hC^n(\cU;\bE_1)$ denotes
all the $\bE_1$-valued,
$n$th continuously differentiable functions on $\cU$, with the usual sup-norm.

\subsection{Fractional calculus}

We begin with a brief review of the deterministic fractional calculus.
We refer to the book
Samko, Kilbas and Marichev \cite{SKM} for an exhaustive survey
on the subject. We first recall some basic definitions.

Let $-\infty<a<b<\infty$, and $\varphi\in L^1([a,b])$. The integrals
%
\begin{eqnarray}
\label{Ia} \bigl(I_{a+}^{\alpha}\varphi\bigr) (x)&=&
\frac{1}{\Gamma(\alpha)}\int_{a}^{x}\frac{\varphi(t)} {
(x-t)^{1-\alpha}}\,\mathrm{d}t,\qquad
  x>a,
\\
\label{Ib} \bigl(I_{b-}^{\alpha}\varphi\bigr) (x)&=&
\frac{1}{\Gamma(\alpha)}\int_{x}^{b}\frac{\varphi(t)} {
(t-x)^{1-\alpha}}\,\mathrm{d}t,\qquad
  x<b,
\end{eqnarray}
are called \textit{fractional integrals of order $\alpha$}, where
$\Gamma
(\cdot)$ is the Gamma-function and $\alpha\in[0,\infty)$.
Both $I^\a_{a+}$ and $I^\a_{b-}$ are the so-called Riemann--Liouville
fractional integrals, and they are often called ``left'' and ``right''
fractional integrals, respectively. We shall denote the image of
$L^p([a,b])$ under the fractional integration operator $I_{a+}^{\a}$
(resp. $I_{b-}^\a$) by $I_{a+}^\a(L^p([a,b]))$ (resp. $I_{b-}^\a
(L^p([a,b]))$). Moreover, in what follows we shall often use
left-fractional integration, which has the following properties:
%
\begin{eqnarray}
\label{Ismg} \bigl[I_{a+}^{\alpha}I_{a+}^{\beta}
\varphi\bigr](\cdot) & = & \bigl[I_{a+}^{\alpha+\beta}\varphi\bigr](
\cdot),
\nonumber
\\[-8pt]\\[-8pt]
t^{\alpha}I^{\beta}_{0+}t^{-\alpha-\beta}I_{0+}^{\alpha}t^{\beta
}
\varphi(\cdot) & = & I_{0+}^{\alpha}I_{0+}^{\beta}
\varphi(\cdot )=I_{0+}^{\alpha+\beta}\varphi(\cdot),\qquad  \alpha>0,
\beta>0.\nonumber
\end{eqnarray}
We note that (\ref{Ismg}) holds for a.e. $x\in[a,b]$. If $\varphi\in
\hC([a,b])$,
then (\ref{Ismg}) holds for all $x\in[a,b]$.
%

The (Riemann--Liouville) fractional derivatives are defined, naturally,
as the\break inverse operator of the
fractional integration. To wit, for any function $f\in L^0([a,b])$, we define
\begin{eqnarray}
\label{Da} \bigl(\mathcal{{D}}_{a+}^{\alpha}f\bigr) (x)&=&
\frac{1}{\Gamma(1-\alpha
)}\frac{\mathrm{d}}{\mathrm{d}x}\int_{a}^{x}
\frac{f(t)}{(x-t)^{\alpha}}\,\mathrm{d}t,
\\
\label{Db} \bigl(\mathcal{{D}}_{b-}^{\alpha}f\bigr) (x)&=&-
\frac{1}{\Gamma(1-\alpha
)}\frac{\mathrm{d}}{\mathrm{d}x}\int_{x}^{b}
\frac{f(t)}{(t-x)^{\alpha}}\,\mathrm{d}t,
\end{eqnarray}
whenever they exist. We call $\cD^\a_{a+}f$ (resp. $\cD^\a_{b-}f$)
the \textit{left (resp. right) fractional derivative of order $\alpha
$}, $0<\alpha<1$. We note that if $f(t)\in\hC^{1}([a,b])$, then it
is easy to verify that (see \cite{SKM}, page 224)
%
\begin{eqnarray}
\label{Da1} \mathcal{{D}}_{a+}^{\alpha}f =\frac{f(x)}{\Gamma(1-\alpha
)(x-a)^{\alpha}}+
\frac{\alpha}{\Gamma(1-\alpha)} \int_{a}^{x}
\frac{f(x)-f(t)}{(x-t)^{1+\alpha}}\,\mathrm{d}t\stackrel{\triangle } {=}D_{a+}^{\alpha}f.
\end{eqnarray}
%
The derivative $D_{a+}^{\alpha}f$ is called Marchaud fractional
derivative. We should note
that the right-hand side of (\ref{Da1}) is not only well-defined for
differentiable
functions, but for example, for function $f(x)$ that is $\beta$-H\"
{o}lder continuous, with
$\beta>\alpha$.
For more general functions,
the fractional Marchaud derivative (\ref{Da1}) should be understood as
(cf. \cite{SKM})
%
\begin{eqnarray}
\label{Da2} D_{a+}^{\alpha}f\stackrel{\triangle} {=}\lim
_{\varepsilon\to
0}D_{a+,\varepsilon
}^{\alpha}f,
\end{eqnarray}
where the limit is in the space $L^{p}$, and
\begin{eqnarray}
\label{Daeps} \bigl[D_{a+,\varepsilon}^{\alpha}f\bigr](x)\stackrel{
\triangle} {=}\frac
{f(x)}{\Gamma(1-\alpha
)(x-a)^{\alpha}}+\frac{\alpha}{\Gamma(1-\alpha)} \int_{a}^{x-\varepsilon}
\frac{f(x)-f(t)}{(x-t)^{1+\alpha}}\,\mathrm{d}t.
\end{eqnarray}

We collect some of the important properties of
the fractional integral and derivative in the follow theorem. The
proofs can be found in
\cite{SKM}.\vadjust{\goodbreak}

\begin{theorem}
\label{fraccalc}
\begin{enumerate}[(iii)]
\item[(i)] For any $\varphi\in L^{1}([a,b])$ and $0<\alpha<1$, it holds that
%
\begin{eqnarray}
\label{DIf} D_{a+}^{\alpha}I_{a+}^{\alpha}
\varphi=\lim_{\varepsilon\rightarrow
0}D_{a+,\varepsilon}^{\alpha}I_{a+}^{\alpha}
\varphi=\mathcal {D}_{a+}^{\alpha}I_{a+}^{\alpha}
\varphi=\varphi.
\end{eqnarray}
%
%

\item[(ii)] For any $f\in I_{a+}^{\alpha}(L^{1}([a,b]))$ and $\a>0$,
it holds that
\begin{eqnarray}
\label{IDf} I_{a+}^{\alpha}D_{a+}^{\alpha}f=I_{a+}^{\alpha}
\mathcal {D}_{a+}^{\alpha}f=f.
\end{eqnarray}
%
\item[(iii)] Let $\psi\in
L^{p}([0,b])$, $b>0$, $1< p<\infty$. Then $\psi$ has the representation
$\psi(x)=I_{0+}^{\alpha}x^{\mu}f(x)$, a.e. $x\in[0,b]$, for some
$f\in L^{p}([0,b])$, $\alpha>0$, and $p(1+\mu)>1$ if and only if
$\psi$ takes one of the following
two forms:
\begin{longlist}[{(b)}]
\item[{(a)}] $\psi(x)=x^{\mu}[I_{0+}^{\alpha}g](x)$, a.e. $x\in
[0,b]$, $g\in L^{p}([0,b])$;\vspace*{2pt}
\item[{(b)}] $\psi(x)=x^{\mu-\varepsilon}[I_{0+}^{\alpha
}x^{\varepsilon}g_{1}](x)$,
a.e. $x\in[0,b]$, $g_{1}\in L^p([0,b])$, $p(1+\varepsilon)>1$.
\end{longlist}
\end{enumerate}
\end{theorem}

\subsection{Fractional Wiener--Poisson space}

We recall that a stochastic process $B^{H}=\{B^{H}_{t}, t\in[0,T]\}$,
defined on a filtered
probability space $(\Omega, \mathcal{F}, \hP; \hF=\{\cF_t\}_{t\ge
0})$, is
called an $\hF$-fractional Brownian motion (fBM) with Hurst
parameter $H\in(0,1)$ if
\begin{longlist}[{(iii)}]
\item[(i)] $B^H$ is a Gaussian process with continuous paths and $B_0^H=0$;

\item[(ii)] for each $t\geq0$, $B^H_{t}$ is $\mathcal{F}_{t}$-measurable and
$\hE B_t^H=0$, for each $t\ge0$;

\item[(iii)] for all $s,t\geq0$, it holds that
%
\begin{eqnarray}
\label{cov} \hE\bigl(B_t^HB_s^H
\bigr)=R_{H}(t,s)=\tfrac{1}{2}\bigl(t^{2H}+s^{2H}-|t-s|^{2H}
\bigr).
\end{eqnarray}
%
\end{longlist}

It follows from (\ref{cov}) that $\hE|B_t^H-B_s^H|^2=|t-s|^{2H}$,
that is,
$B^H$ has stationary increments. Furthermore, by Kolmogorov's
continuity criterion, $B_t^H$ has $\alpha$-H\"older continuous paths
for all $\alpha<H$. In particular, if
$H=1/2$, then $B^H$ becomes a standard Brownian motion; and if $H=1$,
then $\{B_t^1; t\geq0\}$ has the same law as $\{\xi t; t\geq0\}$,
where $\xi$ is an $N(0,1)$ random variable.

In what follows, we shall consider the \textit{canonical space} with
respect to an fBM
or the fractional Wiener space. Let $\Omega^{1}=\hC_0([0,T])$, the
space of all continuous functions, null at zero, and endowed with the
usual sup-norm.
Let $\cF_t^{1}\stackrel{\triangle}{=}\si\{\omega(\cdot\wedge t)|
\omega\in\Omega
^{1}\}$, $t\ge0$, $\cF^{1}\stackrel{\triangle}{=}\cF_T^{1}$, $\hF
^1=\{\cF^{1}_{t},
t\in[0,T]\}$ and $\hP^{B^H}$ is
the probability measure on $(\Omega^{1}, \cF^{1})$ under which the
\textit{canonical process}
\[
B^H_t(\omega)\stackrel{\triangle} {=}\omega(t),\qquad
 (t,\omega )\in[0,T]\times\Omega^{1}
\]
is an fBM of Hurst parameter $H$.



For any $H\in(0,1)$, we define
%
\begin{eqnarray}
\label{RH} R_{H}(t,s)=\int_{0}^{t\wedge s}K_{H}(t,r)K_{H}(s,r)\,\mathrm{d}r,
\end{eqnarray}
where $K_{H}$ is the square integrable kernel given by
%
\begin{eqnarray}
\label{KH} K_{H}(t,s)\stackrel{\triangle} {=}\Gamma\biggl(H+
\frac
{1}{2}\biggr)^{-1}(t-s)^{H-1/2}F\biggl(H-
\frac
{1}{2}, \frac{1}{2}-H, H+\frac{1}{2},1-
\frac{t}{s}\biggr),
\end{eqnarray}
and $F(a,b,c,z)$ is the Gaussian hypergeometric function:
\[
F(a,b,c,z)=\sum_{k=0}^{\infty}
\frac{a^{(k)}b^{(k)}}{c^{(k)}k!}z^{k},  \qquad  a,b\in\hR, |z|<1, c\neq0, -1,\ldots ,
\]
where $a^{(k)}$, $b^{(k)}$, $c^{(k)}$ are the
Pochhammer symbol for the \textit{rising factorial}: $x^{(0)}=1$,
$x^{(k)}=\frac{\Gamma(x+k)}{\Gamma(x)}$.

Now, let $\sE$ be the set of all step functions on $[0,T]$, and let
$\sH$ be the so-called \textit{Reproducing Kernel Hilbert space},
defined as
the closure of $\sE$ with respect to the scalar product
%
\begin{eqnarray}
\label{innerRH} \langle I_{[0,t]},I_{[0,s]}\rangle_{{\sH}}=R_{H}(t,s),\qquad
  s,t\in[0,T].
\end{eqnarray}
%
For any $H\in(0,1)$, we define a linear operator $K_H\dvtx  L^2([0,T])\to
L^2([0,T])$ by
%
\begin{eqnarray}
\label{KHop} [K_{H}f](t)\stackrel{\triangle} {=}\int
_{0}^{t}K_H(t,s)f(s)\,\mathrm{d}s, \qquad   f\in
L^2\bigl([0,T]\bigr),\   t\in[0,T].
\end{eqnarray}
Also, for any $f\in L^0([0,T])$ and $\beta>0$, we shall denote
%
\begin{eqnarray}
\label{lf} [\hspace*{-2pt}[f]\hspace*{-2pt}]^\beta(t)\stackrel{\triangle}
{=}t^\beta f(t),  \qquad  t\in[0,T],
\end{eqnarray}
and
$I_{0+}^{\alpha,\beta}(L^{p}([0,T]))=\{f\in L^0([0,T])\dvtx  [\hspace*{-2pt}[f]\hspace*{-2pt}]^{\beta}\in I_{0+}^{\alpha}(L^{p}([0,T]))\}$.
Then we have the following result (cf., e.g., \cite{DU}, Theorem~2.1, or
\cite{SKM}, Theorem~10.4).

\begin{theorem}
\label{Hprop}
For each $H\in(0,1)$, the operator $K_{H}$ is an isomorphism between
$L^{2}([0,T])$ and $I_{0+}^{H+1/2}(L^{2}([0,T]))$.
Furthermore, it holds that
%
\begin{eqnarray}
\label{KHop1} K_Hf=\lleft\{ %
\begin{array} {l@{\qquad}l}
I_{0+}^{2H}\bigl[\!\!\bigl[I_{0+}^{1/2-H}
[\hspace*{-2pt}[f]\hspace*{-2pt}]^{H-1/2}\bigr]\!\!\bigr]^{1/2-H}, &  H< 1/2,
\\\noalign{\vspace*{2pt}}
I_{0+}^{1}\bigl[\!\!\bigl[I_{0+}^{H-1/2}
[\hspace*{-2pt}[f]\hspace*{-2pt}]^{1/2-H}\bigr]\!\!\bigr]^{H-1/2}, &  H>1/2.
\end{array} %
\rright.
\end{eqnarray}
\end{theorem}

From (\ref{KHop1}) it is easy to check that the inverse operator
$K_{H}^{-1}$ on an absolutely continuous function $h$ satisfies
%
\begin{eqnarray}
\label{KHinv} K_{H}^{-1}h=\lleft\{ %
\begin{array} {l} {\bigl[\!\!\bigl[}I_{0+}^{1/2-H}\bigl[\!\!\bigl[h'\bigr]\!\!\bigr]^{1/2-H}\bigr]\!\!\bigr]^{H-1/2},\\\noalign{\vspace*{2pt}}
\quad  \mbox{if } h'\in L^{1}\bigl([0,T]\bigr), \mbox{and }
H<1/2,
\\\noalign{\vspace*{2pt}}
{\bigl[\!\!\bigl[}D_{0+}^{H-1/2}\bigl[\!\!\bigl[h'
\bigr]\!\!\bigr]^{1/2-H}\bigr]\!\!\bigr]^{H-1/2}, \\\noalign{\vspace*{2pt}}
\quad  \mbox{if }
h'\in I_{0+}^{H-1/2,1/2-H}\bigl(L^{1}
\bigl([0,T]\bigr)\bigr)\cap 
L^{1}\bigl([0,T]\bigr),
\mbox{and } H> 1/2, \end{array} %
\rright.
\end{eqnarray}
where $h^{\prime}$ is the derivative of $h$ (cf., e.g., \cite{SKM}, Theorem~10.6, and \cite{NO}).

Next, let $K_{H}^{*}$ be the adjoint of $K_H$ on $L^2([0,T])$, that is,
for any $f\in\sE,g\in L^2([0,T])$,
\[
\int_0^T \bigl[K_H^{*}f
\bigr](t)g(t)\,\mathrm{d}t=\int_0^T f(t)[K_Hg](t)\,\mathrm{d}t.
\]
Then, it can be shown by Fubini and integration by parts that for any
$f\in\sE$,
\[
\bigl[K^{*}_{H}f\bigr](t)=K_{H}(T,t)\varphi(t)+
\int_{t}^{T}\bigl(f(s)-f(t)\bigr)
\frac
{\partial K_H}{\partial s}(s,t)\,\mathrm{d}s, \qquad   t\in[0,T].
\]
In particular, for $\varphi, \psi\in\sE$, we have (see, e.g., \cite{AMN})
\[
\bigl\langle K_{H}^{*}\varphi, K^{*}_{H}
\psi\bigr\rangle _{L^{2}((0,T))}=\langle\varphi,\psi\rangle_{\mathcal{H}}.
\]
Consequently, the operator $K_{H}^{*}$ is an isometry between the
Hilbert spaces ${\sH}$
and $L^{2}([0,T])$. Furthermore, it can be shown that the process $W=\{
W_{t}, t\in[0,T]\}$ defined by
%
\begin{eqnarray}
\label{W} W_{t}=B^{H}\bigl(\bigl(K_{H}^{*}
\bigr)^{-1}(I_{[0,t]})\bigr)
\end{eqnarray}
is a Wiener process, and the process $B^{H}$ has an integral
representation of the form
%
\begin{eqnarray}
\label{BHrep} B_{t}^{H}=\int_{0}^{t}K_{H}(t,s)\,\mathrm{d}W_{s},\qquad
  t\in[0,T].
\end{eqnarray}

We now turn our attention to the Poisson part. We first consider a
Poisson random
measure $N(\cdot,\cdot)$ on $[0,T]\times\hR$, defined on a given
probability space $(\Omega, \cF, \hP)$, with mean measure
$\hat N(\mathrm{d}t,\mathrm{d}x)=\mathrm{d}t\n(\mathrm{d}x)$, where $\n$ is the L\'evy measure, a $\si
$-finite measure on $\hR^*\stackrel{\triangle}{=}\hR\setminus\{0\}
$ satisfying the
standard integrability condition:
\[
\int_{\hR^*}\bigl(1\wedge|x|^2\bigr)\nu(\mathrm{d}x)<+
\infty.
\]
%
In this paper, we shall be interested in a \textit{Poisson point process}
of class (QL), namely a point process whose counting measure, defined by
$ N_L((0,t]\times A)=\#\{ s\in(0,t]\dvtx \Delta L_s\in A\}=
\sum_{0< s\le t}\mathbf{1}_{\{\Delta L_s\in A\}}$, $t\ge0$, $A\in
\sB
(\hR^*)$,
has a deterministic and continuous compensator
(cf. \cite{IW}). In light of the representation theorem \cite{IW}, Theorem
II-7.4, we shall assume without loss of generality that the process
$L$ takes the following form:
\begin{eqnarray}
\label{L} L_t=\int_{0}^t\int
_{\hR^*}f(s,x)N(\mathrm{d}s,\mathrm{d}x),\qquad    t\ge0,
\end{eqnarray}
where $f\in L^1(\mathrm{d}t\times \mathrm{d}\n)$ is a deterministic function. Then, the
counting measure $N_L(\mathrm{d}t,\mathrm{d}x)$ can be written as
%
\begin{eqnarray}
\label{NL} N_L\bigl((0,t]\times A\bigr)=\int_0^t
\int_{\hR^*}\mathbf{1}_{A}\bigl(f(s,x)
\bigr)N(\mathrm{d}s,\mathrm{d}x), 
\end{eqnarray}
and its compensator is therefore $\h N_L(\mathrm{d}t,\mathrm{d}x)=\hE N_L(\mathrm{d}t,\mathrm{d}x)=f(t,x)\,\mathrm{d}t
\nu(\mathrm{d}x)$. Clearly, if $f(s,x)\equiv g(x)$, then
$L$ is a stationary Poisson point process. In particular, if we assume
that $g(x)\equiv x$ and $\n(\mathrm{d}x)=\l F(\mathrm{d}x)$, where $F(\cdot)$ is a probability
measure on $\hR$, then $L$ is a \textit{compound Poisson process}
with jump intensity $\l$ and jump size distribution $F$.

Throughout this paper, we shall assume that
\begin{eqnarray}
\label{Ui} 
\hE \biggl\{\int
_0^T |L|_t^2\,\mathrm{d}t
+\mathrm{e}^{\beta|\tilde L|_T} \biggr\}<\infty,\qquad    \forall\beta>0, 
\end{eqnarray}
where $|L|_t\stackrel{\triangle}{=}\sum_{0\le s\le t}|\Delta L_s|$
and $|\tilde
L|_t\stackrel{\triangle}{=}\sum_{0\le s\le t}(|\Delta L_s|\vee1)$,
$t\in[0,T]$.

\begin{rem}
\label{remUi}
We note that (\ref{Ui}) contains in particular the compound Poisson
case. Indeed, if $L_t=\sum_{i=1}^{N_t}U_i$, where $N$ is a standard
Poisson process with intensity $\l>0$, and $\{U_i\}$ are
i.i.d. random variables with finite moment generating function
$M_{|U_1|}(t)\stackrel{\triangle}{=}\hE\{\mathrm{e}^{t|U_1|}\}<\infty$,
$\forall t\ge0$.
Then we can easily calculate that
%
\begin{eqnarray}
\label{CPUi} &&\hE \biggl\{\int_0^T
|L|_t^2\,\mathrm{d}t +\mathrm{e}^{\beta|\tilde L|_T} \biggr\}
\nonumber
\\
&&\quad =\frac{(\lambda\hE|U_1|)^{2}T^{3}}{3}+\frac{\lambda\hE\{
|U_1|^2\}T^{2}}{2} + 
\sum
_{k=0}^\infty\hE \bigl\{\mathrm{e}^{\beta\sum_{i=1}^{k} (|U_{i}|\vee
1)}
|N_T=k \bigr\}\frac{(\lambda T)^{k}}{k!}\mathrm{e}^{-\lambda T}
\\
&&\quad =\frac{(\lambda\hE|U_1|)^{2}T^{3}}{3}+\frac{\lambda\hE\{
|U_1|^2\}T^{2}}{2}+ 
\mathrm{e}^{\lambda T(\hE[\mathrm{e}^{\beta( |U_{1}|\vee1)}]-1)}<\infty.\nonumber
\end{eqnarray}
%
\end{rem}

We can also consider the canonical space for a given Poisson point
process of class (QL). Let $\Omega^{2}=\hD([0,T])$, the
space of all real-valued, c\`adl\`ag (right-continuous with left limit)
functions, endowed
with the Skorohod topology, and let $\cF_t^{2}\stackrel{\triangle
}{=}\si\{\omega
(\cdot\wedge t)| \omega\in\Omega^{2}\}$, $t\ge0$, $\cF
^{2}\stackrel{\triangle}{=}\cF
_T^{2}$, $\hF^2=\{\cF^{2}_{t}, t\in[0,T]\}$.
Let $\hP^{L}$ be the law of the process $L$ on $\hD([0,T])$. Then,
the coordinate process,
by a slight abuse of notations,
\[
L_t(\omega)=\omega(t), \qquad (t,\omega)\in[0,T]\times
\Omega^{2},
\]
is a Poisson point process, defined on $(\Omega^{2 },\cF^{2},\hP^{L})$,
whose compensated counting measure is
$\h N_L(\mathrm{d}t,\mathrm{d}z)=\hE[N_L(\mathrm{d}t,\mathrm{d}z)]=f(t,z)\n(\mathrm{d}z)\,\mathrm{d}t$, where $\n$ is a L\'
evy measure and (\ref{Ui}) holds.

Combining the discussions above, we now consider two canonical spaces
$(\Omega^1,\cF^1,\hP^{B^H};\hF^1)$ and
$(\Omega^2,\cF^2,\hP^{L};\hF^2)$, where $\Omega^1=\hC([0,T])$
and $\Omega
^2=\hD([0,T])$. We define
the \textit{fractional Wiener--Poisson space} to simply be the product space:
%
\begin{eqnarray}
\label{WP} \Omega&\stackrel{\triangle} {=}&\Omega^1\times
\Omega^2; \qquad \cF \stackrel{\triangle} {=}\cF^1\otimes
\cF^2; \nonumber\\[-8pt]\\[-8pt]
 \hP &\stackrel{\triangle} {=}&\hP^{B^H}\otimes
\hP^{L};\qquad   \cF _t\stackrel{\triangle} {=}
\cF^1_t\otimes\cF ^2_t,\qquad  t\in[0,T].\nonumber
\end{eqnarray}
We write the element of $\Omega$ as $\omega=(\omega^1,\omega^2)\in
\Omega$. Then, the two marginal coordinate processes defined by
%
\begin{eqnarray}
\label{coord} B^H_t(\omega)\stackrel{\triangle} {=}
\omega^1(t), \qquad   L_t(\omega )\stackrel{\triangle} {=}
\omega ^2(t), \qquad  (t,\omega)\times[0,T]\times\Omega,
\end{eqnarray}
will be the fractional Brownian motion and Poisson point process,
respectively, with
the given laws. Note that under our assumptions $B^H$ and
$L$ are always independent (cf., e.g., \cite{IW}, Theorem II-6.3).
Also, we can assume without loss of generality that
the filtration $\hF$ is right continuous, and is augmented by all the
$\hP$-null
sets so that it
satisfies the \textit{usual hypotheses}.

To end this section, we recall that if $\cX$ is a metric space, $X$ is
a $\cX$-valued Gaussian random variable, and $g(\cdot)$ is a
seminorm on $\cX$, such that
and $\hP(g(X)<\infty)>0$. Then it follows from the Fernique Theorem
(cf. \cite{Fernique}) that there exists $\varepsilon>0$ such that $
\hE[ \exp(\l g^2(X))] <\infty$, for all $0<\lambda<\varepsilon$.
It is then easy to see that for all $0<\rho<2$, one has
%
\begin{eqnarray}
\label{Fernique} \hE\bigl[ \exp\bigl(\lambda g^\rho(X)\bigr)\bigr] <
\infty, \qquad  \forall\l>0.
\end{eqnarray}
This fact is useful in our analysis, similar to, for example, \cite{NO}.


%

\section{The problem}\label{sec3}

In this paper, we are interested in the following stochastic
differential equation with additive noise:
%
\begin{eqnarray}
\label{SDE1} X_{t}=x+\int_{0}^{t}b(s,X_{s})\,\mathrm{d}s
+B_{t}^{H}-L_{t},\qquad  t\in[0,T],
\end{eqnarray}
where $b$ is a Borel function on $[0,T]\times\hR$, $B^H$ is an fBM
with Hurst parameter
$H\in(0,1)$ and $L$ is a Poisson point process of class (QL), both
defined on some filtered probability space
$(\Omega, \cF, \hP; \hF)$. We assume that $B^H$ and $L$ are both
$\hF
$-adapted, and they are independent.
We often consider the filtration generated by $(B^H, L)$, denoted by
$\hF^{(B^H,L)}=\{\cF^{(B^H,L)}_t\dvtx t\ge0\}$ where
%
\begin{eqnarray}
\label{FBL} \cF^{(B^H,L)}_t\stackrel{\triangle} {=}\si\bigl
\{\bigl(B^H_s,L_s\bigr)\dvtx  0\le s\le t\bigr
\}, \qquad  t\ge0,
\end{eqnarray}
and we assume that $\hF^{(B^H,L)}$ is augmented by all the $\hP$-null
sets so that it satisfies the \textit{usual hypotheses}. As usual, we
have the following definitions of solutions to the SDE (\ref{SDE1}).

\begin{defn}
\label{strongdef}
Let $(\Omega, \cF, \hP)$ be a complete probability space on which are
defined an fBM $B^H$,
$H\in(0,1)$, and a Poisson point process $L$, independent of $B^H$ and
of class (QL). A process $X$ defined on $(\Omega,\cF, \hP)$ is
called a
strong solution to (\ref{SDE1}) if
\begin{enumerate}[(ii)]
\item[(i)] $X$ is $\hF^{(B^H,L)}$-adapted;

\item[(ii)] $X$ satisfies (\ref{SDE1}), $\hP$-almost surely.
\end{enumerate}%
\end{defn}

\begin{defn}
\label{weakdef}
A seven-tuple $(\Omega, \mathcal{F}, P, \hF, X, B^{H}, L)$ is called a
weak solution to (\ref{SDE1}) if
\begin{enumerate}[(iii)]
\item[(i)] $(\Omega, \mathcal{F}, P; \hF)$ is a filtered probability space;

\item[(ii)] $B^{H}$ is an $\hF$-fBM, and $L$ is an $\hF$-Poisson point
process of class (QL);

\item[(iii)] $(X, B^{H}, L)$ satisfies (\ref{SDE1}), $\hP$-almost surely.
\end{enumerate}%
\end{defn}

For simplicity, we often say that $(X, B^H, L)$ (or simply $X$) is a
weak solution to
(\ref{SDE1}) without specifying the associated probability space
$(\Omega
, \cF, \hP;\hF)$ when the context is clear. It is readily seen from
(\ref{SDE1}) that if $(X, B^H, L)$ is a weak solution, then $\hF
^{(B^H, L)}\subseteq\hF^X$. The well-known example of
Tanaka indicates that the converse is not necessarily true, even in the
case when
$H=1/2$ and $L\equiv0$.

Throughout this paper, we shall make use of the following \textit
{standing assumptions}:

\begin{assum}
\label{Assump1} The function $b\dvtx  [0,T]\times\hR\mapsto\hR$
satisfies the following assumptions for $H\in(0,1/2)$ and $H\in
(1/2,1)$, respectively:
\begin{enumerate}[(ii)]
\item[(i)] If $H< 1/2$, then for some $0<\rho\le1$ and $K>0$, it holds that
%
\begin{eqnarray}
\label{Hle12} \bigl|b(t,x)\bigr|\leq K\bigl(1+|x|^{\rho}\bigr),\qquad
\forall(t,x)\in[0,T]\times\hR.
\end{eqnarray}

\item[(ii)] If $H>1/2$, then $b$ is H\"{o}lder-$\gamma$ continuous in $t$
and H\"older-$\a$ in
$x$, where $\gamma>H-1/2$, and $1-\frac{1}{2H}<\alpha<1$.
That is, for some $K>0$,
%
\begin{eqnarray}
\label{Hge12} \bigl|b(t,x)-b(s,y)\bigr|\leq K\bigl(|x-y|^{\alpha}+|t-s|^{\gamma}
\bigr),  \qquad \forall (t,x), (s,y)\in[0,T]\times\hR.
\end{eqnarray}
\end{enumerate}%
\end{assum}

\begin{rem}
\label{rem1}
(1) We note that in the case when $H<1/2$ we do not require any
regularity on the coefficient $b$. To discuss the well-posedness under
such a weak condition on
the coefficient, is only possible due to the presence of the ``noises''
$B^H$ and
$L$ (see also \cite{NO} for the case when $L\equiv0$), and it is
quite different from the theory of ordinary differential equations, for example.

(2) Compared to \cite{NO}, we require that $b$ grows only sub-linearly
in the case $H<1/2$. This is
due to the possible infinite jumps of $L$. In fact, Remark~\ref
{rhoge12} below shows that the problem could
be ill-posed if $\rho>1/2$. Such a constraint can be removed when $L$
has only finitely many jumps.
\end{rem}

We end this section by making the following observation. Denote $\tilde
X=X+L$, and
\[
\tilde b(t,x,\omega)\stackrel{\triangle} {=}b\bigl(t,x-L_t(\omega)
\bigr), \qquad  (t,x,\omega)\in [0,T]\times\hR\times\Omega.
\]
Then the SDE (\ref{SDE1}) becomes
%
\begin{eqnarray}
\label{SDE2} \tilde X_{t}=x+\int_{0}^{t}
\tilde b(s,\tilde X_{s})\,\mathrm{d}s +B_{t}^{H},\qquad    t
\in[0,T].
\end{eqnarray}
Thus the problem is reduced to the case studied by \cite{NO}, except
that the coefficient $\tilde b$ is now random. However, if we consider
the problem
on the canonical Wiener--Poisson space in which $(B^H_t(\omega
),L_t(\omega))=(\omega^1(t), \omega^2(t))$, $t\in[0,T]$, then we
can formally consider the SDE (\ref{SDE2}) as one
on $(\Omega^1, \cF^1, \hP^{B^H})$:
%
\begin{eqnarray}
\label{SDE3} \tilde X_{t}=x+\int_{0}^{t}
b^{\omega^2}(s,\tilde X_{s})\,\mathrm{d}s +B_{t}^{H},\qquad
  t\in[0,T],
\end{eqnarray}
where $b^{\omega^2}(t,x)\stackrel{\triangle}{=}b(t,x-\omega
^2(t))=\tilde
b(t,x,\omega^2)$, for each fixed $\omega^2\in\Omega^2$. In other
words, we can apply the result of \cite{NO} to obtain the
well-posedness for each $\omega^2\in\Omega^2$, provided that the
coefficient $b^{\omega^2}$ satisfies the assumptions in \cite{NO}.
We should note, however, that such a seemingly simple argument is
actually rather difficult to implement, especially for the weak
solution case, due to some subtle measurability issues caused by the
lack of
regularity of $b$ in the case $H<1/2$, and the discontinuity of the
paths of $L$
(whence $\tilde b$ in the temporal variable $t$),
in the case $H>1/2$.

\section{Existence of a weak solution ($H<1/2$)}\label{sec4}

In this section, we shall validate the argument presented at the end of
the last section, in the case $H<1/2$. Namely, we shall prove that the
SDE (\ref{SDE2}) possesses a weak solution, along the lines of the
arguments of
\cite{NO}.

Recall from Assumption~\ref{Assump1} that in the case $H<1/2$ the
function $b$
satisfies (\ref{Hle12}).
Consider the canonical Wiener--Poisson space $(\Omega,\cF,\hP)$, where
$\hP=\hP^{B^H}\otimes\hP^L$, with a given Hurst parameter $H\in
(0,1/2)$, a L\'evy measure $\nu(\mathrm{d}z)$,
and a deterministic function $f\dvtx [0,T]\times\hR\mapsto\hR$ so that
$\h N_L(\mathrm{d}t,\mathrm{d}z)=\hE[N_L(\mathrm{d}t,\mathrm{d}z)]=f(t,z)\n(\mathrm{d}z)\,\mathrm{d}t$ satisfies
(\ref{Ui}). Let $(B^{H}, L)$ be the canonical process. Define
$u_t\stackrel{\triangle}{=}-b(t, B_t^H-L_t+x)$
and
%
\begin{eqnarray}
\label{v} v_{t}\stackrel{\triangle} {=}-K_{H}^{-1}
\biggl(\int_{0}^{\cdot}b\bigl(r,
B_{r}^{H}-L_{r}+x\bigr)\,\mathrm{d}r \biggr)
(t)=K_{H}^{-1} \biggl(\int_{0}^{\cdot
}u_r\,\mathrm{d}r
\biggr) (t), \qquad   t\in[0,T],
\end{eqnarray}
where $K_H^{-1}$ is defined by (\ref{KHinv}). We have the following lemma.

\begin{lem}
\label{Girsanov}
Assume $H<1/2$ and (\ref{Hle12}) is in force with $0<\rho<1/2$. Then
the process $v$ defined by (\ref{v}) enjoys the following properties:
\begin{enumerate}[(2)]
\item[(1)] $\hP\{v\in L^2([0,T])\}=1$;

\item[(2)] $v$ satisfies the Novikov condition:
%
\begin{eqnarray}
\label{Novikov} \hE \biggl\{\exp \biggl(\frac{1}{2}\int
_{0}^{T}|v_t|^{2}\,\mathrm{d}t
\biggr) \biggr\} <\infty.
\end{eqnarray}
\end{enumerate}%
Furthermore, if $L$ has only finitely many jumps, then the results hold
under (\ref{Hle12}) for any $\rho\in(0,1]$.
\end{lem}

\begin{pf} (1) In what follows, we denote $C>0$ to be a generic
constant depending only on the coefficient $b$, the constants in
Assumption~\ref{Assump1}, and the Hurst parameter $H$; and is allowed
to vary from line to line. Since $H<1/2$, and (\ref{Hle12}) holds,
some simple computation, together with assumption (\ref{Ui}), shows that
\begin{eqnarray*}
\hE\int_0^T|u_t|^2
\,\mathrm{d}t&=&\hE\int_{0}^{T}\bigl|b\bigl(t,B_{t}^{H}-L_{t}+x
\bigr)\bigr|^2\,\mathrm{d}s\leq C\hE\int_{0}^{T}
\bigl(1+\bigl|B_t^{H}-L_t+x\bigr|\bigr)^{2}\,\mathrm{d}t
\\
&\leq& C \biggl[ \bigl(1+|x|\bigr)^{2}T+ \hE\int_{0}^{T}\bigl|B_t^{H}\bigr|^{2}\,\mathrm{d}t+
\hE \int_{0}^{T}|L|_t^{2}\,\mathrm{d}t
\biggr]
\\
&=& C \biggl[\bigl(1+|x|\bigr)^{2}T+\frac{T^{2H+1}}{2H+1} +\hE\int
_0^T|L|^2_t\,\mathrm{d}t \biggr] <
\infty.
\end{eqnarray*}
%
Therefore, $\int_{0}^{T}|u_t|^{2}\,\mathrm{d}s<\infty$, $\hP$-a.s. Since
$H<1/2$, $[\hspace*{-2pt}[u]\hspace*{-2pt}]^{1/2-H}$ belongs
to $L^2([0,T])$, $\hP$-a.s. as well. Thus, applying \cite{SKM}, Theorem~5.3,
$I_{0+}^{1/2-H}[\hspace*{-2pt}[u]\hspace*{-2pt}]^{1/2-H}\in L^q([0,T])$, $\hP$-a.s., for
some $q=\frac{2}{1-2(1/2-H)}=\frac{1}{H}>2$. In particular,\vspace{2pt}
$I_{0+}^{1/2-H}[\hspace*{-2pt}[u]\hspace*{-2pt}]^{1/2-H}\in L^2([0,T])$, $\hP$-a.s. Let
$N\subset\Omega$ be the
exceptional $\hP$-null set. Then for any $\omega\notin N$, we can
apply Theorem~\ref{fraccalc}(iii)(a)
to find $h^\omega\in L^{2}([0,T])$ such that
\[
\bigl[I_{0+}^{1/2-H}[\hspace*{-2pt}[u]\hspace*{-2pt}]^{1/2-H}(\omega
)\bigr](t)=t^{1/2-H}\bigl[I_{0+}^{1/2-H}h^\omega
\bigr](t),  \qquad \omega\notin N.
\]
Now recall from (\ref{KHinv}) we see that this implies that for each
$\omega\notin N$, it holds that
\[
K_{H}^{-1} \biggl(\int_{0}^{\cdot}u_{r}(
\omega)\,\mathrm{d}r \biggr)=I_{0+}^{1/2-H}h^\omega.
\]
Thus, applying \cite{SKM}, Theorem~5.3, again we have $K_{H}^{-1}(\int_{0}^{\cdot}u_{r}(\cdot)\,\mathrm{d}r)\in L^q([0,T])$, $\hP$-a.s., for some
$q=\frac{2}{1-2(1/2-H)}=\frac{1}{H}>2$. In
particular, (1) holds.\vspace*{1pt}

(2) Using the Assumption~\ref{Hle12} again we have, $\hP$-almost surely,
\begin{eqnarray*}
|v_{s}|&=&\bigl|s^{H-1/2}I_{0+}^{1/2-H}[\hspace*{-2pt}[u]\hspace*{-2pt}]^{1/2-H}(s)\bigr|
\\
&=&C s^{H-1/2} \biggl|\int_{0}^{s}(s-r)^{-1/2-H}r^{1/2-H}b
\bigl(r,B_{r}^{H}-L_{r}+x\bigr)\,\mathrm{d}r \biggr|
\\
&\leq& C T^{1/2-H}\bigl(1+|x|^{\rho}+\bigl\|B^H
\bigr\|_{\infty}^{\rho
}+|L|_{T}^{\rho}\bigr),
\end{eqnarray*}
where $\|B^{H}\|_{\infty}\stackrel{\triangle}{=}\sup_{0\leq s \leq
T}|B_{s}^{H}|$.
Note that $L$ and $B^H$ are independent we have
%
\begin{eqnarray}
\label{Novikov0} &&\hE \biggl\{\exp \biggl(\frac{1}2\int
_0^T|v_{t}|^2\,\mathrm{d}t \biggr)
\biggr\}\nonumber\\[-8pt]\\[-8pt]
&&\quad \le \mathrm{e}^{CT^{2-2H}(1+|x|^{2\rho})}\hE \bigl\{\exp \bigl(CT^{2-2H}
\bigl\|B^H\bigr\| _{\infty}^{2\rho} \bigr) \bigr\} \hE\bigl
\{\mathrm{e}^{CT^{2-2H}|L|_{T}^{2\rho}}\bigr\}.\nonumber
\end{eqnarray}
Note that $2\rho<1$ by (\ref{Hle12}) in Assumption~\ref{Assump1}, we have
\begin{eqnarray}
\label{Novikov1} \hE\bigl\{\mathrm{e}^{CT^{2-2H}|L|_{T}^{2\rho}} \bigr\}\leq\hE\bigl\{ \mathrm{e}^{CT^{2-2H}(|L|_{T}+1)}
\bigr\}<\infty,
\end{eqnarray}
thanks to (\ref{Ui}).
Note that $\rho<1/2 $ also guarantees that $\hE \{\exp
(CT^{2-2H}\|B^H\|_{\infty}^{2\rho} ) \}<\infty$
for all $T>0$ with $\cX=\hC([0,T])$, $X=B^H$, and $g(\cdot)=\|
\cdot\|_\infty$ in (\ref{Fernique}).
This, together with (\ref{Novikov0}) and (\ref{Novikov1}), proves
(\ref{Novikov}).

Finally, note that if $L$ has
only finitely many jumps, then $\Delta L_t=0$ for all but finitely many
$t\in[0,T]$. Thus (\ref{Novikov1}) holds for
all $\rho\in(0,1]$. This proof is now complete.
\end{pf}
%

\begin{rem}
\label{rhoge12}
We note that unlike the finite jump case (see also \cite{NO} for the
continuous case) where we only
assume $0<\rho\le1$, in general it is necessary to assume $\rho<1/2$
to guarantee the finiteness of $\hE\{\mathrm{e}^{|L|^{2\rho}_T}\}$. In fact,
if $\rho>1/2$, then even in the
simplest standard Poisson case $L_t\equiv N_t$ we have
\[
\hE \mathrm{e}^{(N_T)^{2\rho}}=\sum_{n=0}^{\infty}\mathrm{e}^{n^{2\rho}}
\frac
{\lambda^{n}}{n !}\mathrm{e}^{-\lambda}. 
\]
If we denote $a_{n}=\mathrm{e}^{n^{2\rho}}\frac{\lambda^{n}}{n !}$, then $\ln
a_{n}=n^{2\rho}+n\ln\lambda-\ln n!$. Since $\ln n!<n \ln n$, and
\[
\lim_{n\rightarrow\infty}\frac{n\ln n }{n^{2\rho}+n \ln\lambda
}=0,
\]
 a simple calculation
then shows that
\begin{eqnarray*}
\lim_{n\to\infty}\ln a_n&=&\lim_{n\to\infty}
\bigl\{n^{2\rho}+n \ln \lambda-\ln n!\bigr\}
\\
&=&\lim_{n\to\infty}\bigl\{n^{2\rho}+n \ln\lambda\bigr\}
\biggl\{1-\frac{\ln
n! }{n^{2\rho}+n \ln\lambda} \biggr\}= +\infty.
\end{eqnarray*}
That is, $a_{n}\rightarrow+\infty$, and consequently $\hE
\mathrm{e}^{(N_T)^{2\rho}}=\infty$.
\end{rem}


We can now construct a weak solution to (\ref{SDE1}), in the case
$H<1/2$, as follows. Define
%
\begin{eqnarray}
\label{tidBH} \tilde{B}_{t}^{H}\stackrel{\triangle}
{=}B_{t}^{H}-\int_{0}^{t}b
\bigl(s, B_{s}^{H}-L_{s}+x\bigr)\,\mathrm{d}s
=B_{t}^{H}+\int_{0}^{t}u_{s}\,\mathrm{d}s,\qquad
  t\in[0,T].
\end{eqnarray}
Using the representation (\ref{BHrep}), we can write
\begin{eqnarray*}
\tilde{B}_{t}^{H}=B_{t}^{H}+\int
_{0}^{t}u_{s}\,\mathrm{d}s=\int
_{0}^{t}K_{H}(t,s)\,\mathrm{d}W_{s}+
\int_{0}^{t}u_{s}\,\mathrm{d}s=\int
_{0}^{t}K_{H}(t,s)\,\mathrm{d}
\tilde{W}_{s},
\end{eqnarray*}
where
%
\begin{eqnarray}
\label{WienerInt} \tilde{W}_{t}=W_{t}+\int
_{0}^{t}\biggl(K_{H}^{-1}
\biggl(\int_{0}^{.}u_{s}\,\mathrm{d}s\biggr) (r)
\biggr)\,\mathrm{d}r=W_{t}+\int_{0}^{t}v_r\,\mathrm{d}r.
\end{eqnarray}
By Lemma~\ref{Girsanov}, the process $v$
satisfies the Novikov condition (\ref{Novikov}). Thus, if we define a
new probability
measure $\tilde P$ on the canonical fractional Wiener--Poisson space
$(\Omega, \cF)$ by
%
\begin{eqnarray}
\label{tildeP} \frac{\mathrm{d}\tilde P}{\mathrm{d}P}\stackrel{\triangle} {=}\exp \biggl\{-\int
_{0}^{T}v_{s}\,\mathrm{d}W_{s}-
\frac
{1}{2}\int_{0}^{T}v_{s}^{2}\,\mathrm{d}s
\biggr\},
\end{eqnarray}
%
then, under $\tilde{P}$, $\tilde W$ is an $\hF$-Brownian motion, and
$\tilde{B}^{H}$ is an $\hF$-fractional Brownian motion with Hurst
parameter H (cf. Decreusefond and \"Ustunel \cite{DU}).

Furthermore, since $B^H$ and $L$ are independent, we can easily check,
by following the
arguments of Brownian case (cf., e.g., \cite{S}, Theorem~124, \cite
{IW}, Theorem II-6.3) that $L_{t}$ is still a Poisson point process
of class (QL) with same parameters, and is
independent of $\tilde{B}^{H}$.
We now define $X_t=x+B^{H}_t-L_t$, $t\in[0,T]$. Then, it follows from
(\ref{tidBH})
that
%
\begin{eqnarray}
\label{weakX} \tilde B^H_t=(X_t-x+L_t)-
\int_0^tb(t,X_s)\,\mathrm{d}s, \qquad   t
\in[0,T].
\end{eqnarray}
In other words, $(\Omega,\cF,\tilde\hP,\hF, X, \tilde B^H, L)$ is a
weak solution of (\ref{SDE1}).
%
%
That is, we have proved the following theorem.

\begin{theorem}
\label{weakHle12}
Assume $H<1/2$ and that the assumptions of Lemma~\ref{Girsanov} are in
force. Then
for any $T>0$, the SDE (\ref{SDE1}) has at least one weak solution on $[0,T]$.
\end{theorem}


%
%

\section{Existence of a weak solution ($H>1/2$)}\label{sec5}

In this section, we study the existence of the weak solution in the
case when $H>1/2$. We note that even though the coefficient $b$
is H\"older continuous in both variables by Assumption~\ref
{Assump1}(ii) (\ref{Hge12}), the coefficient $\tilde b$
of the reduced SDE (\ref{SDE2}) will have discontinuity on the
variable $t$, thus the Assumption~\ref{Assump1}(ii) is no longer
valid for
$\tilde b$, and therefore the results of \cite{NO} cannot be
applied directly. We shall, however, using the
same scheme as in the last section to prove the existence of the weak
solution, although the arguments is much more involved.

We begin with some preparations. Let $(\Omega, \mathcal{F}, \hP, \hF
)$ be the canonical fractional Wiener--Poisson space, and let $(B^{H},
L)$ be the canonical process. For fixed $x\in\hR$, consider again the
process
\begin{eqnarray*}
u_t(\omega)=-b\bigl(t, B_t^{H}(
\omega)-L_t(\omega)+x\bigr)=-b\bigl(t, \omega ^1(t)-
\omega^2(t)+x\bigr),  \qquad (t,\omega)\in[0,T]\times\Omega,
\end{eqnarray*}
and define $v_t(\omega)=K_{H}^{-1} (\int_{0}^{\cdot
}u_{r}(\omega)\,\mathrm{d}r )(t)$, $(t,\omega)\in[0,T]\times\Omega$, where
$K_H^{-1}$ is given by (\ref{KHinv}) in the case $H>1/2$. As in the
previous section, we shall
again argue that Lemma~\ref{Girsanov} holds. The main difference
between our case and \cite{NO}, however, is that the paths of
$u$ are discontinuous despite the Assumption~\ref{Assump1}(ii), thus
the fractional calculus will need to be modified.

We first note that, by the Fubini theorem,
\[
\hP\bigl\{v\in L^2\bigl([0,T]\bigr)\bigr\}=\int_{\Omega^2}
\hP^{B^H} \biggl\{\int_0^T\bigl|v_s
\bigl(\omega^1,\omega^2\bigr)\bigr|^2\,\mathrm{d}s<\infty
\biggr\} \hP^L\bigl(\mathrm{d}\omega^2\bigr).
\]
Thus to show $\hP\{v\in L^2([0,T])\}=1$, it suffices to show that, for
$\hP^L$-a.e., $\omega^2\in\Omega^2$, it holds that
\[
\hP^{B^H} \biggl\{\int_0^T\bigl|v^{\omega^2}_s
\bigl(\omega^1\bigr)\bigr|^2\,\mathrm{d}s<\infty \biggr\}=1,
\]
where $v^{\omega^2}_s(\omega^1)\stackrel{\triangle}{=}v_s(\omega
^1,\omega^2)$ is
the ``$\omega^2$-section'' of $v_t$. But in light of (\ref{KHinv}),
we need first show that, for $\hP^L$-a.e. $\omega^2\in\Omega^2$,
$u^{\omega^2}\in I_{0+}^{H-1/2, 1/2-H}(L^{1}([0, T]))\cap
L^1([0,T])$, $\hP^{B^H}$-a.s., where
%
\begin{eqnarray}
\label{uw2} u^{\omega^2}_t(\omega_{1})\stackrel{
\triangle} {=}u_t\bigl(\omega _{1},\omega ^2
\bigr)=-b^{\omega^2,x}\bigl(t,B^H_t(
\omega_{1})\bigr),\qquad  \bigl(t,\omega^1\bigr)\in [0,T]
\times\Omega^1
\end{eqnarray}
and
%
\begin{eqnarray}
\label{bOx} b^{\omega^2,x}(t,y)\stackrel{\triangle} {=}b\bigl(t, y-
\omega^2(t)+x\bigr), \qquad  (t,y)\in [0,T]\times\hR.
\end{eqnarray}

Since we are considering only the canonical process $L(\omega
)=L(\omega^2)=\omega^2$, which is a Poisson process under $\hP^L$
and thus
does not have fixed time jumps (i.e., $\hP^L\{\Delta L_t \neq0\}=0$,
$\forall t\ge0$). We can, modulo a $\hP^L$-null set, assume without
of generality that
$\omega^2$ is piecewise constant, and
jumps at $0<\sigma_{1}(\omega^2)< \cdots<\sigma_{N_T(\omega
^2)}(\omega^2)<T$, where $N_t(\omega^2)$ denotes the number of
jumps of $L(\omega^2)$ up
to time $t>0$. For notational convenience in what follows, we shall
also denote
$\sigma_{0}(\omega^2)=0$, $\si_{N_T(\omega^2)+1}(\omega^2)=T$,
although they do not represent jump times. 
Then by Assumption~\ref{Assump1}(ii) we see that $t\mapsto b^{\omega
^2,x}(t,B_t^H)$ is $\mu$-H\"older continuous
on every interval $(\sigma_{i}, \sigma_{i+1})$, $i=0,1,\ldots
,N_T(\omega^2)$,
with $\mu=H-\frac{1}{2}+\varepsilon$ for some $\varepsilon>0$.
Thus, by virtue of Theorem~6.5 in~\cite{SKM}, $u^{\omega^2}\in
I_{\sigma_{i}+}^{H-1/2}(L^{2}(\si_{i}, \si_{i+1}))$, $\hP
^{B^H}$-a.s., for all $i=0,\ldots,N_T(\omega^2)$. It then follows from
Theorem~13.11 of \cite{SKM} that $u^{\omega^2}\in
I_{0+}^{H-1/2}(L^{2}([0, T]))$,
$\hP^{B^H}$-a.s.
Therefore, there exists a $\hP^{B^H}$-null set $N\subset\Omega^1$, so
that for any $\omega^1\notin N$, we can apply Theorem~\ref{fraccalc}(iii)(a)
or Lemma~3.2 in \cite{SKM} to find a function $h^{\omega^1,\omega
^2}\in L^{2}([0, T])$, such that:
\[
\bigl[\!\!\bigl[u^{\omega^2}\bigr]\!\!\bigr]^{1/2-H}\bigl(t,
\omega^1\bigr)=t^{1/2-H}u^{\omega
^2}_t\bigl(
\omega^1\bigr)=I_{0+}^{H-1/2}t^{1/2-H}h^{\omega^1,\omega
^2}(t),
 \qquad  t\in[0,T].
\]
That is,
$u^{\omega^2}\in I_{0+}^{H-1/2, 1/2-H}(L^{1}([0, T]))$, $\hP
^{B^H}$-a.s. On the other hand, since
$u^{\omega^2}\in I_{0+}^{H-1/2}(L^{2}[0,\allowbreak  T])$ implies $u^{\omega
^2}\in L^{2}([0, T])$, thanks to Theorem~5.3
of \cite{SKM}, we conclude that
(\ref{KHinv}) holds with $h(\cdot)=\int_{0}^{\cdot}u_{r}\,\mathrm{d}r$, $\hP
^{B^H}$-a.s. That is,
$v_{t}= K_{H}^{-1}(\int_{0}^{\cdot}u_{r}\,\mathrm{d}r)(t)$, $t\in[0,T]$, belongs\vspace*{1.5pt}
to $L^{2}([0,T])$, $\hP^{B^H}$-a.s.
Note that the argument is valid for $\hP^L$-a.e. $\omega^2\in\Omega
^2$, we obtain that $\hP\{v\in L^2([0,T])\}=1$.
%
We now prove an analogue of Lemma~\ref{Girsanov} for the case $H>1/2$.

\begin{lem}
\label{Girsanov1}
Assume that $H>1/2$, and that Assumption~\ref{Assump1}\textup{(ii)} holds with
$1-\frac{1}{2H}<\a< 1-H$.
Then the conclusion of Lemma~\ref{Girsanov} remains valid.

Furthermore, if $L$ has only finitely many jumps, then the constraint
$\a<1-H$ can be removed.
\end{lem}

\begin{pf} We have already argued that the process $v_{t}=
K_{H}^{-1}(\int_{0}^{\cdot}u_{r}\,\mathrm{d}r)(t)$, $t\in[0,T]$, satisfies
$\hP\{v\in L^2([0,T])\}=1$ in the beginning of this section. We shall
show that the process $v$ also satisfies the Novikov condition (\ref
{Novikov}), whence part (2) of Lemma~\ref{Girsanov}.

To this end, first note that on the
canonical space $\Omega^2=\hD([0,T])$, and under the probability
$\hP
^L$, the canonical
process $L(\omega)=\omega^2$ is a Poisson point process of class (QL).
Now, for fixed $T>0$, denote $\Omega^2_n\stackrel{\triangle}{=}\{
\omega^2\dvtx  N_T(\omega
^2)=n\}$ for $n=0,1,\ldots $;
and for $\omega^2\in\Omega^2_n$, again denote
$0<\sigma_{1}(\omega^2)<\cdots<\si_n(\omega^2)<T$ be the
jump times of $L(\omega^2)$, and $\si_0(\omega^2)=0$, $\si
_{n+1}(\omega^2)=T$. Finally, denote $S_k(\omega^2)\stackrel
{\triangle}{=}\sum_{i=1}^k\Delta L_{\si_i}(\omega^2)$, $k=1,2,\ldots, $ and
$S_0(\omega^2)=0$. In what follows,\vspace{2pt} we often suppress the variable
$\omega^2$ when the context
is clear.

Now recall from (\ref{KHinv}) that, for $H>1/2$,
%
\begin{eqnarray}
\label{KHinv-u} v^{\omega^2}_t=K_{H}^{-1}
\biggl(\int_{0}^{\cdot}u^{\omega
^2}_{r}\,\mathrm{d}r
\biggr) (t)= t^{H-1/2}D_{0+}^{H-1/2}\bigl[\!\!
\bigl[u^{\omega^2}\bigr]\!\!\bigr]^{1/2-H}(t), \qquad   t\in[0,T].
\end{eqnarray}
We shall calculate $D_{0+}^{H-1/2}[\hspace*{-2pt}[u^{\omega^2}]\hspace*{-2pt}]^{1/2-H}$ for
$\omega^2\in\Omega^2_n$, for each
$n=0, 1,2,\ldots\, $. To see this, fix $n\in\hN$, and let $\omega
^2\in
\Omega^2_n$. For notational simplicity,
in what follows we denote
%
\begin{eqnarray}
\label{uw2k} u^{\omega^2,k}_t\bigl(\omega^1\bigr)=-b
\bigl(t, B^H_t\bigl(\omega^1
\bigr)-S_{k-1}\bigl(\omega ^2\bigr)+x\bigr),\qquad  \bigl(t,
\omega^1\bigr)\in[0,T]\times\Omega^1, k\ge1,
\end{eqnarray}
so that
$ u^{\omega^2}_t=\sum_{k=1}^{n+1} u^{\omega^2,k}_t\mathbf
{1}_{[\si_{k-1}(\omega^2),\si_k(\omega^2))}(t)$, $t\in[0,T]$,
$\hP^1\mbox{-a.s.}$
Then, for\vspace{2pt}
$t\in[0,\si_1(\omega^2))$, by definition (\ref{Da2}) and (\ref
{Daeps}) with $p=2$ we
have
%
\begin{eqnarray}
\label{DHv0} &&D_{0+}^{H-1/2}\bigl[\!\!\bigl[u^{\omega^2}
\bigr]\!\!\bigr]^{1/2-H}(t)
\nonumber
\\
&&\quad =\frac{1}{\Gamma(3/2-H)}\frac{[\hspace*{-2pt}[u^{\omega^2,1}]\hspace*{-2pt}]^{1/2-H}(t)}{t^{H-1/2}}\\
&&\qquad {}+\frac{H-1/2}{\Gamma(3/2-H)} \int
_{0}^{t}\frac{[\hspace*{-2pt}[u^{\omega^2,1}]\hspace*{-2pt}]^{1/2-H}(t)-[\hspace*{-2pt}[u^{\omega^2,1}]\hspace*{-2pt}]^{1/2-H}(r)}{(t-r)^{H+1/2}}\,\mathrm{d}r\nonumber
\\
&&\quad \stackrel{\triangle} {=}\Phi_1(t).\nonumber
\end{eqnarray}
Similarly, for $\sigma_{k-1}(\omega^2)\le t<\sigma_{k}(\omega^2)$
with $1<k\le n+1$, we have
%
\begin{eqnarray}
\label{DHvi} &&D_{0+}^{H-1/2}\bigl[\!\!\bigl[u^{\omega^2}
\bigr]\!\!\bigr]^{1/2-H}(t)
\nonumber
\\
&&\quad =\frac{1}{\Gamma(3/2-H)}\frac{[\hspace*{-2pt}[u^{\omega^2}]\hspace*{-2pt}]^{1/2-H}(t)}{t^{H-1/2}}\nonumber\\
&&\qquad {}+\frac{H-1/2}{\Gamma(3/2-H)} \int
_{0}^{t}\frac{[\hspace*{-2pt}[u^{\omega^2}]\hspace*{-2pt}]^{1/2-H}(t)-[\hspace*{-2pt}[u^{\omega
^2}\bigl]\!\!\bigl]^{1/2-H}(r)}{(t-r)^{H+1/2}}\,\mathrm{d}r
\nonumber
\\[-8pt]\\[-8pt]
&&\quad =\frac{1}{\Gamma(3/2-H)}\frac{[\hspace*{-2pt}[u^{\omega^2,k}]\hspace*{-2pt}]^{1/2-H}(t)}{t^{H-1/2}}
\nonumber\\
&&\qquad {}+\frac{H-1/2}{\Gamma(3/2-H)}\sum_{i=1}^{k-1} \int
_{\si_{i-1}}^{\sigma_{i}}\frac{[\hspace*{-2pt}[u^{\omega^2,k}]\hspace*{-2pt}]^{1/2-H}(t)-
[\hspace*{-2pt}[u^{\omega^2,i}]\hspace*{-2pt}]^{1/2-H}(r)}{(t-r)^{H+1/2}}\,\mathrm{d}r
\nonumber
\\
&&\qquad {}+\frac{H-1/2}{\Gamma(3/2-H)} \int_{\sigma_{k-1}}^{t}
\frac{[\hspace*{-2pt}[u^{\omega^2,k}]\hspace*{-2pt}]^{1/2-H}(t)-[\hspace*{-2pt}[u^{\omega^2,k}]\hspace*{-2pt}]^{1/2-H}(r)}{(t-r)^{H+1/2}}\,\mathrm{d}r \stackrel{\triangle} {=}\Phi_k(t).\nonumber 
\end{eqnarray}
Consequently, we obtain the following formula:
%
\begin{eqnarray}
\label{DHv} D_{0+}^{H-1/2}\bigl[\!\!\bigl[u^{\omega^2}
\bigr]\!\!\bigr]^{1/2-H}(t)=\sum_{k=1}^{n+1}
\Phi_{k}(t) \b1_{[\si_{k-1}(\omega^2),\si_{k}(\omega^2))}(t), 
 \qquad  t
\in[0,T), \hP^1\mbox{-a.s.}
\end{eqnarray}
%
That is,
%
\begin{eqnarray}
\label{vo2} v^{\omega^2}_{t}=t^{H-1/2}D_{0+}^{H-1/2}
\bigl[\!\!\bigl[u^{\omega^2}\bigr]\!\!\bigr]^{1/2-H}(t)=t^{H-1/2}
\sum_{i=1}^{n+1}\Phi_k(t)
\b1_{[\sigma_{k-1}(\omega^2)< t\leq\sigma_{k}(\omega^2))}(t),
\end{eqnarray}
%
where $\Phi_k$'s are defined by (\ref{DHv0}) and (\ref{DHvi}). We
now estimate each term in (\ref{vo2}).
Note that for $t\in[\sigma_{k-1}, \sigma_{k})$ we have
\begin{eqnarray*}
&&\frac{H-\sfrac{1}2}{\Gamma(\sfrac{3}{2}-H)}\sum_{i=1}^{k-1} \int
_{\si_{i-1}}^{\sigma_{i}}\frac{[\hspace*{-2pt}[u^{\omega^2,k}]\hspace*{-2pt}]^{1/2-H}(t)}{(t-r)^{H+1/2}}\,\mathrm{d}r\\
&&\quad = \frac{1}{\Gamma(\sfrac{3}2-H)}
\biggl\{\frac{1}{(t-\sigma
_{k-1})^{H-\sfrac{1}{2}}}- \frac{1}{t^{H-\sfrac{1}{2}}} \biggr\}\bigl[\!\!
\bigl[u^{\omega^2,k}\bigr]\!\!\bigr]^{1/2-H}(t). 
\end{eqnarray*}
It then follows from (\ref{DHvi}) that, for $t\in[\si_{k-1},\si_k)$,
\begin{eqnarray*}
t^{H-1/2}\Phi_k(t)& =& t^{H-1/2} \biggl\{
\frac{1}{\Gamma(3/2-H)}\frac
{[\hspace*{-2pt}[u^{\omega^2,k}]\hspace*{-2pt}]^{1/2-H}(t)}{t^{H-1/2}}
\\
&&\hphantom{t^{H-1/2} \{}{}+\frac{H-1/2}{\Gamma(3/2-H)}\sum_{i=1}^{k-1} \int
_{\si_{i-1}}^{\sigma_{i}}\frac{[\hspace*{-2pt}[u^{\omega^2,k}]\hspace*{-2pt}]^{1/2-H}(t)-
[\hspace*{-2pt}[u^{\omega^2,i}]\hspace*{-2pt}]^{1/2-H}(r)}{(t-r)^{H+1/2}}\,\mathrm{d}r
\nonumber
\\
&&\hphantom{t^{H-1/2} \{}{}+\frac{H-1/2}{\Gamma(3/2-H)} \int_{\sigma_{k-1}}^{t}
\frac{[\hspace*{-2pt}[u^{\omega^2,k}]\hspace*{-2pt}]^{1/2-H}(t)-[\hspace*{-2pt}[u^{\omega^2,k}]\hspace*{-2pt}]^{1/2-H}(r)}{(t-r)^{H+1/2}}\,\mathrm{d}r \biggr\}
\nonumber
\\
& =& C_1^H\frac{t^{H-1/2}[\hspace*{-2pt}[u^{\omega^2,k}]\hspace*{-2pt}]^{1/2-H}(t)}{(t-\sigma_{k-1})^{H-1/2}}-C_2^Ht^{H-1/2}
\sum_{i=1}^{k-1}\int_{\sigma_{i-1}}^{\sigma_{i}}
\frac{[\hspace*{-2pt}[u^{\omega
^2,i}]\hspace*{-2pt}]^{1/2-H}(t)} {
(t-r)^{H+1/2}}\,\mathrm{d}r
\\
&&{}+C_2^Ht^{H-1/2}\sum
_{i=1}^{k-1}\int_{\sigma_{i-1}}^{\sigma
_{i}}
\frac{[\hspace*{-2pt}[u^{\omega^2,i}]\hspace*{-2pt}]^{1/2-H}(t)
-[\hspace*{-2pt}[u^{\omega^2,i}]\hspace*{-2pt}]^{1/2-H}(r)}{(t-r)^{H+1/2}}\,\mathrm{d}r
\\
&&{}+C_2^Ht^{H-1/2}\int_{\sigma_{k-1}}^{t}
\frac{[\hspace*{-2pt}[u^{\omega
^2,k}]\hspace*{-2pt}]^{1/2-H}(t)-[\hspace*{-2pt}[u^{\omega^2,k}]\hspace*{-2pt}]^{1/2-H}(r)} {
(t-r)^{H+1/2}}\,\mathrm{d}r\stackrel{\triangle} {=}A^{k}(t)+B^{k}(t),
\end{eqnarray*}
where $C_1^H\stackrel{\triangle}{=}\frac{1}{\Gamma(3/2-H)}$,
$C_2^H\stackrel{\triangle}{=}\frac
{H-1/2}{\Gamma(3/2-H)}=(H-1/2)C_1^H$, and
%
\begin{eqnarray}
\label{Ak} A^{k}(t)&\stackrel{\triangle} {=}&C_1^H
\frac{t^{H-1/2}[\hspace*{-2pt}[u^{\omega
^2,k}]\hspace*{-2pt}]^{1/2-H}(t)}{(t-\sigma_{k-1})^{H-1/2}}- C_2^Ht^{H-1/2}\sum
_{i=1}^{k-1}\int_{\sigma_{i-1}}^{\sigma_{i}}
\frac
{[\hspace*{-2pt}[u^{\omega^2,i}]\hspace*{-2pt}]^{1/2-H}(t)} {
(t-r)^{H+1/2}}\,\mathrm{d}r,
\\
\label{Bi} B^{k}(t)&\stackrel{\triangle} {=}&C_2^Ht^{H-1/2}
\sum_{i=1}^{k-1}\int_{\sigma
_{i-1}}^{\sigma_{i}}
\frac{[\hspace*{-2pt}[u^{\omega^2,i}]\hspace*{-2pt}]^{1/2-H}(t)
-[\hspace*{-2pt}[u^{\omega^2,i}]\hspace*{-2pt}]^{1/2-H}(r)}{(t-r)^{H+1/2}}\,\mathrm{d}r
\nonumber
\\[-8pt]\\[-8pt]
&&{}+C_2^Ht^{H-1/2}\int_{\sigma_{k-1}}^{t}
\frac{[\hspace*{-2pt}[u^{\omega
^2,k}]\hspace*{-2pt}]^{1/2-H}(t)-[\hspace*{-2pt}[u^{\omega^2,k}]\hspace*{-2pt}]^{1/2-H}(r)} {
(t-r)^{H+1/2}}\,\mathrm{d}r.\nonumber
\end{eqnarray}
It is readily seen that (suppressing $\omega=(\omega^1,\omega^2)$'s)
\begin{eqnarray}\label{Ai}
\bigl|A^k(t)\bigr| &=& \Biggl|C_1^H\sum
_{i=1}^{k-1}b\bigl(t,B_{t}^{H}-S_{i-1}+x
\bigr) \biggl[\frac{1}{(t-\sigma_{i-1})^{H-1/2}}-\frac{1} {
(t-\sigma_{i})^{H-1/2}} \biggr]
\nonumber
\\
&&\hphantom{\Biggl|}{}+C_1^H\frac{b(t,B^H_{t}-S_{k-1}+x)}{(t-\sigma_{k-1})^{H-1/2}} \Biggr|
\nonumber
\\
&\leq& \Biggl|C_1^H\sum_{i=1}^{k-1}
\bigl[b\bigl(t,B_{t}^{H}-S_{i-1}+x\bigr)-b
\bigl(t,B_{t}^{H}+x\bigr)\bigr] \biggl[\frac{1}{(t-\sigma_{i-1})^{H-1/2}}-
\frac{1} {
(t-\sigma_{i})^{H-1/2}} \biggr]
\nonumber
\\
&&\hphantom{\Biggl|}{}+C_1^H\frac{(b(t,B^H_{t}-S_{k-1}+x)-b(t,B^H_{t}+x)}{(t-\sigma
_{k-1})^{H-1/2}} \Biggr|
\nonumber
\\
&&\hphantom{\Biggl|}{}+ \Biggl|C_1^H\sum_{i=1}^{k-1}b
\bigl(t,B_{t}^{H}+x\bigr) \biggl[\frac{1}{(t-\sigma
_{i-1})^{H-1/2}}-
\frac{1} {
(t-\sigma_{i})^{H-1/2}} \biggr]\nonumber\\[-8pt]\\[-8pt]
&&\hphantom{\Biggl|{}+\Biggl|}{}+C_1^H\frac{b(t,B^H_{t}+x)}{(t-\sigma
_{k-1})^{H-1/2}} \Biggr|\nonumber
\\
 &\leq& C_1^H\max_{1\leq i \leq
k}\bigl|b
\bigl(t,B_{t}^{H}-S_{i-1}+x\bigr)-b
\bigl(t,B_{t}^{H}+x\bigr)\bigr|
\nonumber
\\
&&{} \times\Biggl|\sum_{i=1}^{k-1} \biggl[
\frac{1}{(t-\sigma_{i})^{H-1/2}}-\frac{1} {
(t-\sigma_{i-1})^{H-1/2}} \biggr]+\frac{1}{(t-\sigma
_{k-1})^{H-1/2}}
\Biggr|\nonumber\\
&&{}+C_1^H t^{1/2-H} \bigl|b\bigl(t,B_{t}^{H}+x
\bigr) \bigr|
\nonumber
\\
&\leq&C(t-\sigma_{k-1})^{1/2-H}|L|_{T}^{\alpha
}+Ct^{1/2-H}
\bigl(\bigl|b(0,x)\bigr|+|t|^{\gamma}+\bigl\|B^{H}\bigr\|_{\infty}^{\alpha}
\bigr),\nonumber
\end{eqnarray}
where
$C$ is a generic constant depending on $H$, $\a$, and $K$, thanks to
Assumption~\ref{Assump1}. On the other hand, we write $
B^{k}(t)=-C(B_{1}^{k}(t)+B_{2}^{k}(t))$, where
\begin{eqnarray*}
B_{1}^{k}(t)&=&t^{H-1/2}\sum
_{i=1}^{k-1} \biggl[b\bigl(t,B_{t}^{H}-S_{i-1}+x
\bigr)\int_{\sigma_{i-1}}^{\sigma_{i}}\frac{t^{1/2-H}
-r^{1/2-H}}{(t-r)^{1/2+H}}\,\mathrm{d}r
\\
&&\hphantom{t^{H-1/2}\sum
_{i=1}^{k-1} \biggl[}{}+ \int_{\sigma_{i-1}}^{\sigma_{i}}\frac{b(t,
B_{t}^{H}-S_{i-1}+x)-b(r,
B_{t}^{H}-S_{i-1}+x)}{(t-r)^{1/2+H}}r^{1/2-H}\,\mathrm{d}r
\biggr]
\\
&&{}+t^{H-1/2}b\bigl(t,B_{t}^{H}-S_{k-1}+x
\bigr)\int_{\sigma_{k-1}}^{t}\frac
{t^{1/2-H}-r^{1/2-H}}{(t-r)^{1/2+H}}\,\mathrm{d}r
\\
&&{}+ t^{H-1/2}\int_{\sigma_{k-1}}^{t}
\frac{b(t, B_{t}^{H}-S_{k-1}+x)-b(r,
B_{t}^{H}-S_{k-1}+x)}{(t-r)^{1/2+H}}r^{1/2-H}\,\mathrm{d}r ,
\end{eqnarray*}
and
\begin{eqnarray*}
B_{2}^{k}(t)&=&t^{H-1/2}\sum
_{i=1}^{k-1}\int_{\sigma_{i-1}}^{\sigma
_{i}}
\frac{b(r,B_{t}^{H}-S_{i-1}+x)-b(r,
B_{r}^{H}-S_{i-1}+x)}{(t-r)^{1/2+H}}r^{1/2-H}\,\mathrm{d}r
\\
&&{}+t^{H-1/2}\int_{\sigma_{k-1}}^{t}
\frac
{b(r,B_{t}^{H}-S_{k-1}+x)-b(r, B_{r}^{H}-S_{k-1}+x)}{(t-r)^{1/2+H}}r^{1/2-H}\,\mathrm{d}r.
\end{eqnarray*}
Then, it is easy to see that, for each fixed $0<\varepsilon<H-\frac
{H-1/2}{\alpha}$ (recall Assumption~\ref{Assump1}(ii)), and denoting
$G\stackrel{\triangle}{=}{\sup}_{0\leq t<r\leq T} \frac
{|B_{t}^{H}-B_{r}^{H}|}{|t-r|^{H-\varepsilon}}$,
we have
%
\begin{eqnarray}
\label{Bk2} \bigl|B_{2}^{k}(t)\bigr|&\leq&t^{H-1/2}\sum
_{i=1}^{k-1}\int_{\sigma
_{i-1}}^{\sigma_{i}}
\frac{|B_{t}^{H}-B_{r}^{H}|^{\alpha}} {
(t-r)^{1/2+H}}r^{1/2-H}\,\mathrm{d}r\nonumber\\
&&{}+t^{H-1/2}\int_{\sigma_{k-1}}^{t}
\frac
{|B_{t}^{H}-B_{r}^{H}|^{\alpha}}{(t-r)^{1/2+H}}r^{1/2-H}\,\mathrm{d}r
\\
&=&t^{H-1/2}\int_{0}^{t}
\frac{|B_{t}^{H}-B_{r}^{H}|^{\alpha
}}{(t-r)^{1/2+H}}r^{1/2-H}\,\mathrm{d}r\leq Ct^{1/2-H+\alpha(H-\varepsilon
)}G^{\alpha}.\nonumber
\end{eqnarray}
Furthermore, by the same argument as in (\ref{Ai}) we also have
\begin{eqnarray}\label{Bk1}
\bigl|B_{1}^{k}(t)\bigr|&=&t^{H-1/2} \max
_{1\leq i\leq
k}\bigl|b\bigl(t,B_{t}^{H}-S_{i-1}+x
\bigr)\bigr|
\nonumber
\\
&&{}\times \Biggl[\sum_{i=1}^{k-1}\int
_{\sigma_{i-1}}^{\sigma_{i}}\frac
{r^{1/2-H}-t^{1/2-H}}{(t-r)^{1/2+H}}\,\mathrm{d}r +\int
_{\sigma_{k-1}}^{t}\frac
{r^{1/2-H}-t^{1/2-H}}{(t-r)^{1/2+H}}\,\mathrm{d}r \Biggr]
\nonumber
\\
 &&{}+K t^{H-1/2} \Biggl[\sum_{i=1}^{k-1}
\int_{\sigma_{i-1}}^{\sigma
_{i}}\frac{|t-r|^{\gamma}}{(t-r)^{1/2+H}}r^{1/2-H}\,\mathrm{d}r+
\int_{\sigma_{k-1}}^{t}\frac{|t-r|^{\gamma
}}{(t-r)^{1/2+H}}r^{1/2-H}\,\mathrm{d}r
\Biggr]
\nonumber
\\
&\leq&\bigl[\bigl|b(0,x)\bigr|+K\bigl(|t|^{\gamma}+\bigl|B_{t}^{H}\bigr|^{\alpha}+|L_{T}|^{\alpha
}
\bigr)\bigr]t^{H-1/2}\int_{0}^{t}
\frac
{r^{1/2-H}-t^{1/2-H}}{(t-r)^{1/2+H}}\,\mathrm{d}r
\\
&&{}+Kt^{H-1/2}\int_{0}^{t}
\frac{|t-r|^{\gamma
}}{(t-r)^{1/2+H}}r^{1/2-H}\,\mathrm{d}r
\nonumber
\\
&\leq& C\bigl\{\bigl[\bigl|b(0,x)\bigr|+|t|^{\gamma}+\bigl|B_{t}^{H}\bigr|^{\alpha
}+|L_{T}|^{\alpha}
\bigr]t^{1/2-H}+t^{\gamma+1/2-H}\bigr\}
\nonumber
\\
&\leq&Ct^{1/2-H}\bigl[\bigl|b(0,x)\bigr|+|t|^{\gamma}+\bigl\|B^{H}
\bigr\|_{\infty}^{\alpha
}+|L|_{T}^{\alpha}\bigr].\nonumber
\end{eqnarray}
Combining (\ref{Bk2}) and (\ref{Bk1}), we have for any $t\in[0,T]$,
%
\begin{eqnarray}
\label{Bk} \bigl|B^k(t)\bigr|\le Ct^{1/2-H}\bigl[\bigl|b(0,x)\bigr|+|t|^{\gamma}+
\bigl\|B^{H}\bigr\|_{\infty
}^{\alpha}+|L|_{T}^{\alpha}+t^{\a(H-\varepsilon)}G^{\a}
\bigr].
\end{eqnarray}

Now, combining (\ref{Ai}) and (\ref{Bk}), and denoting $\hE_n[\cdot
]=\hE[\cdot|N_T=n]$, we have
\begin{eqnarray}
\label{Novikov2a} &&\hE \biggl\{\exp \biggl\{\frac{1}{2} \int
_{0}^{T}v^{2}(t)\,\mathrm{d}t \biggr\} \biggr\}
\nonumber
\\
&&\quad = \sum_{n=0}^{\infty}\hE_n \Biggl\{\exp
\Biggl\{\frac{1}{2}\sum_{k=1}^{n}
\int_{\sigma_{k-1}}^{\sigma_{k}} t^{2H-1}\Phi
^{2}_k(t)\,\mathrm{d}t\\
&&\qquad \hphantom{\sum_{n=0}^{\infty}\hE_n \Biggl\{\exp
\Biggl\{}{}+\frac{1}{2}\int_{\sigma_{n}}^{T}
t^{2H-1}\Phi ^{2}_{n+1}(t)\,\mathrm{d}t \Biggr
\}\Biggr
\}P(N_T=n)\nonumber
\\
&&\quad =  \sum_{n=0}^{\infty}\hE_{n} \Biggl\{\exp
\Biggl\{C \sum_{k=1}^{n+1}\int
_{\sigma_{k-1}}^{\sigma_{k}}\bigl(|A^{k}(t)|+
|B^k(t)|\bigr)^{2}\,\mathrm{d}t \Biggr\}\Biggr
\}P(N_T=n).\nonumber
\end{eqnarray}
%
By (\ref{Ai}) and (\ref{Bk}) and using the fact\vspace*{-2pt}
\[
\frac{\sum_{i=1}^{n+1}x_{i}^{2-2H}}{n+1}\leq \biggl(\frac{\sum_{i=1}^{n+1} x_{i}}{n+1} \biggr)^{2-2H},\qquad
x_{i}>0,
\]
we have\vspace*{-2pt}
%
\begin{eqnarray}
\label{Ak+Bk} &&\sum_{k=1}^{n+1}\int
_{\sigma_{k-1}}^{\sigma_{k}}\bigl(\bigl|A^{k}(t)\bigr|+
\bigl|B^k(t)\bigr|\bigr)^{2}\,\mathrm{d}t\nonumber\\
&&\quad \le 
C\sum
_{k=1}^{n+1}\int_{\sigma_{k-1}}^{\sigma_{k}}(t-
\sigma _{k-1})^{1-2H}|L|_{T}^{2\alpha} \,\mathrm{d}t
\nonumber
\\
&& \qquad {}   + C\int_{0}^{T}t^{1-2H}
\bigl(\bigl|b^{2}(0,x)\bigr|+|t|^{2\gamma}+\bigl\|B^{H}
\bigr\|_{\infty
}^{2\alpha} +t^{2\alpha(H-\varepsilon)}G^{2\alpha}\bigr)\,\mathrm{d}t
\nonumber\\[-8pt]\\[-8pt]
&&\quad \leq 
C\sum_{k=1}^{n+1}(
\sigma_{k}-\sigma_{k-1})^{2-2H}|L|_{T}^{2\alpha}\nonumber\\
&&\qquad {}+C\int_{0}^{T}t^{1-2H}
\bigl(\bigl|b^{2}(0,x)\bigr| +|t|^{2\gamma}+\bigl\|B^{H}
\bigr\|_{\infty}^{2\alpha}+t^{2\alpha
(H-\varepsilon)}G^{2\alpha}\bigr)\,\mathrm{d}t
\nonumber
\\
&&\quad \leq 
C(n+1)^{2H-1}|L|_{T}^{2\alpha}+C
\bigl[1+\bigl\|B^{H}\bigr\|_{\infty}^{2\alpha
}+G^{2\alpha}
\bigr].\nonumber
\end{eqnarray}
Putting (\ref{Ak+Bk}) into (\ref{Novikov2a}), we obtain\vspace*{-2pt}
%
\begin{eqnarray}
\label{Novikov2} &&\hE \bigl\{\mathrm{e}^{\sfrac{1}{2} \int_{0}^{T}v^{2}(t)\,\mathrm{d}t} \bigr\}\nonumber\\[-8pt]\\[-8pt]
&&\quad \leq\hE \bigl\{\exp\bigl\{C
\bigl[1+\bigl\|B^{H}\bigr\|_{\infty}^{2\alpha}+G^{2\alpha}
\bigr]\bigr\} \bigr\}\hE \bigl\{\exp\bigl\{C (N_{T}+1)^{2H-1}|L|_{T}^{2\alpha}
\bigr\} \bigr\}.\nonumber
\end{eqnarray}
By the same argument as Lemma~\ref{Girsanov}, it is easy to prove that
$\hE \{\mathrm{e}^{C\|B^{H}\|_{\infty}^{2\alpha}+G^{2\alpha}} \}
<\infty$.

We need to show that $\hE \{\mathrm{e}^{C (N_{T}+1)^{2H-1}|L|_{T}^{2\alpha
}} \}<\infty$. Note that $\alpha<1-H$ in Assumption~\ref
{Assump1}(ii) implies that $2H-1+2\alpha<1$, and recall $\tilde L$
from (\ref{Ui}), we have\vspace*{-2pt}
%
\begin{eqnarray}
\label{est2} \hE\exp\bigl\{C (N_{T}+1)^{2H-1}|L|_{T}^{2\alpha}
\bigr\}&\leq& \hE\exp \Biggl\{C \Biggl(\sum_{i=1}^{N_T}\bigl(|
\Delta L_{\si_i}|\vee1\bigr)+1 \Biggr)^{2H-1+2\alpha} \Biggr\}
\nonumber
\\[-8pt]\\[-8pt]
&\leq& \hE\exp \Biggl\{C \Biggl(\sum_{i=1}^{N_T}|
\Delta\tilde L_{\si
_i}|+1 \Biggr) \Biggr\}<\infty. \nonumber
\end{eqnarray}
Therefore, we can show that $\hE \{\mathrm{e}^{\sfrac{1}{2} \int
_{0}^{T}v^{2}(t)\,\mathrm{d}t} \}<\infty$.

Finally, note that if $L$ has only finitely many jumps, then $ |\Delta
L_{\si_i}|=0$ for all but finitely many $i$'s. Thus, (\ref{est2})
always holds for any $\a>0$. The proof is complete.
\end{pf}
%

\begin{rem}
We observe that $1-\frac{1}{2H}<\a<1-H$ implies $H<\frac{\sqrt{2}}2$.
This is again due to the presence of possible infinite number of jumps.
We note that a similar constraint $H<\frac{1+\sqrt{5}}4$ was also
placed in \cite{MN}, where only finitely many jumps were considered.
But in that case we need only
$1-\frac{1}{2H}<\a< 1$, thus our result is still much stronger than
that of \cite{MN}.
\end{rem}

We have the following analogues of Theorem~\ref{weakHle12}.

\begin{theorem}
\label{weakHge12}
Assume $H>1/2$ and that the assumptions in Lemma~\ref{Girsanov1} are
in force. Then the SDE (\ref{SDE1}) has at least one weak solution on $[0,T]$.
\end{theorem}

\section{Uniqueness in law and pathwise uniqueness}\label{sec6}

In this section, we study the uniqueness of the weak solution. We shall
first show that the weak solutions to (\ref{SDE1}) are unique in law.
The argument is very similar to that of \cite{NO},
we describe it briefly.

Let $(X, B^{H},L)$ be a weak solutions of (\ref{SDE1}), defined on
some probability
space $(\Omega,\cF,\hP; \hF)$, with the existence interval
$[0,T]$. Let
$W$ be the $\hF$-Brownian motion such that
%
\begin{eqnarray}
B^{H}_{t}=\int_{0}^{t}K_{H}(t,s)\,\mathrm{d}W_{s},\qquad
  t\in[0,T].
\end{eqnarray}
Define
%
\begin{eqnarray}
\label{v-2} v_{t}=K_{H}^{-1} \biggl(\int
_{0}^{\cdot}b(r, X_{r} )\,\mathrm{d}r \biggr) (t),\qquad
  t\in[0,T],
\end{eqnarray}
and let us assume that $v$ satisfies the assumption (1) and (2) in
Lemma~\ref{Girsanov}. Then applying the Girsanov theorem we see that
the process
$\tilde{W}_t =W_{t} +\int_{0}^{t}v_{s} \,\mathrm{d}s$, $t\in[0,T]$,
is an $\hF$-Brownian motion under the new probability measure $\tilde
{\hP}$, defined by
%
\begin{eqnarray}
\frac{\mathrm{d}\tilde{\hP} }{\mathrm{d}\hP}=\xi_{T}(X)\stackrel{\triangle} {=}\exp \biggl\{-
\int_{0}^{T}v_{t} \,\mathrm{d}W_{t} -
\frac{1}{2}\int_{0}^{T}|v_{t}|^{2}\,\mathrm{d}t
\biggr\}.
\end{eqnarray}
Thus $\tilde B^H_t\stackrel{\triangle}{=}\int_{0}^{t}K_{H}(t,s)\,\mathrm{d}\tilde W_{s}$, $t\in
[0,T]$, is an fBM under $\tilde\hP$, and it holds that
\[
X_{t}+L_{t}-x=\int_{0}^{t}b(s,X_{s})\,\mathrm{d}s+B^{H}_{t}=
\int_{0}^{t}K_{H}(t,s)\,\mathrm{d}
\tilde{W}_{s}=\tilde B^H_t,  \qquad  t\in[0,T].
\]
Since under the Girsanov transformation the process $L$ remains a
Poisson point process with the same parameters, and is automatically
independent of the Brownian motion $\tilde W$ under $\tilde\hP$ (cf.
\cite{IW}, Theorem II-6.3), we can then write $X$ as the independent
sum of $\tilde B^H$ and $-L$:
\[
X_t = x+\tilde B^H_t-L_t,\qquad
t\in[0,T].
\]
Since the argument above can be applied to any weak solution, we have
essentially proved the following weak uniqueness result.
%

\begin{theorem}
\label{Uweak1}
Suppose that the assumptions of Lemma~\ref{Girsanov} (resp. Lemma~\ref
{Girsanov1}) for $H<1/2$ (resp. $H>1/2$) are in force. Then two weak
solutions of SDE (\ref{SDE1}) must have the same law, over their
common existence interval $[0,T]$.
\end{theorem}

\begin{pf} We need only to show that the adapted process $v$ defined
by (\ref{v-2}) satisfies
(1) and (2) in Lemma~\ref{Girsanov}. In what follows we
let $C>0$ denote a generic constant depending only on the constants
$H$, $K$, $\a$, $\gamma$ in Assumption~\ref{Assump1} and
$T>0$, and is allowed to vary from line to line. In the case $H<\frac
{1}{2}$, denoting $u=b(\cdot,X_\cdot)$, for any $t\in[0,T]$ we have
\begin{eqnarray*}
\hE\int_{0}^{t}|u_r|^2\,\mathrm{d}r&=&
\hE\int_{0}^{t}\bigl|b(r,X_{r})\bigr|^2\,\mathrm{d}r
\leq C\hE\int_{0}^{t}\bigl(1+|X_{r}
|^{2}\bigr)\,\mathrm{d}r
\\
&\le& C\hE\int_{0}^{t} \biggl[1+|x|^{2}+
\biggl|\int_{0}^{r}b(s,X_{s} )\,\mathrm{d}s
\biggr|^{2}+\bigl|B_{r}^{H}\bigr|^{2}+|L_{r}|^{2}
\biggr]\,\mathrm{d}r
\\
&\le& C \biggl\{\hE\int_{0}^{t}r\int
_{0}^{r}|u_s|^2\,\mathrm{d}s\,\mathrm{d}r+
\bigl(1+|x|^2\bigr)t+\frac{t^{2H+1}}{2H+1}+\hE\int_0^T|L|_T^2\,\mathrm{d}r
\biggr\}
\\
&\le&C_L \biggl\{\bigl(1+|x|^2
\bigr)+\int_{0}^{t}\hE\int_{0}^{r}|u_s|^2\,\mathrm{d}s\,\mathrm{d}r
\biggr\},
\end{eqnarray*}
where $C_L>0$ depends on $C$ and $L$, thanks to (\ref{Ui}). Thus by
Growall's inequality, we obtain
\begin{eqnarray*}
\hE\int_0^T|u_s|^2\,\mathrm{d}s=
\hE\int_{0}^{T}\bigl|b(s,X_{s})\bigr|^2\,\mathrm{d}s
\leq C_L\bigl(1+|x|^2\bigr)\mathrm{e}^{C_LT}<\infty.
\end{eqnarray*}
Then, by the same argument as Lemma~\ref{Girsanov}, we can check that
$v=K_H^{-1}(\int_0^\cdot u_r\,\mathrm{d}r)$ satisfies (1) of Lemma~\ref{Girsanov}.
Furthermore, similarly to the proof Lemma~\ref{Girsanov} we can obtain that
\begin{eqnarray*}
|v_{s}|\leq C_LT^{1/2-H}\bigl(1+\|X
\|_{\infty}^{\rho}\bigr),
\end{eqnarray*}
where $\|X\|_{\infty}\stackrel{\triangle}{=}\sup_{0\leq s \leq
T}|X_s|$. Applying
Grownall's inequality again it is easy to show that
%
\begin{eqnarray}
\label{Xinfinity} \|X \|_{\infty}\leq\bigl(|x|+\bigl\|B^{H}
\bigr\|_{\infty}+C_LT+|L|_{T}\bigr)\mathrm{e}^{C_LT},
\end{eqnarray}
which then leads to (2) of Lemma~\ref{Girsanov}.

We now assume $H>\frac{1}{2}$. Following the same argument of Lemma~\ref{Girsanov1}, it suffices
to show that between two jump times of $L$, the process $u = b(\cdot
,X_\cdot)\in I_{{\si_{k-1}}+}^{H-1/2}(L^{2}([\si_{k-1},\si_k)))$,
$\hP$-almost surely. But note that between two jumps we have, by
Assumption~\ref{Assump1},
\begin{eqnarray*}
&&\bigl|b(t,X_{t} )-b(s,X_{s} )\bigr|\\
&&\quad \leq C\bigl
\{|t-s|^{\gamma}+|X_{t} -X_{s} |^{\alpha}\bigr\}
\\
&&\quad \leq C\biggl\{ |t-s|^{\gamma}+ \biggl|\int_{s}^{t}b(u,X_{u}
)\,\mathrm{d}u \biggr|^{\alpha}+\bigl|B_{t}^{H}-B_{s}^{H}\bigr|^{\alpha}
\biggr\}
\\
&&\quad \leq C\biggl\{ |t-s|^{\gamma}+\biggl |\int_{s}^{t}
\bigl(\bigl|b(0,x)\bigr|+|u|^{\gamma
}+|X_{u}-x|^{\alpha}\bigr)\,\mathrm{d}u
\biggr|^{\alpha}+\bigl|B_{t}^{H}-B_{s}^{H}\bigr|^{\alpha
}
\biggr\}
\\
&&\quad \leq C\bigl\{|t-s|^{\gamma}+\bigl(\bigl|b(0,x)\bigr|+|T|^{\gamma}+\|X
\|_{\infty
}^{\alpha}+|x|^\a\bigr)|t-s|^{\alpha}+\bigl|B_{t}^{H}-B_{s}^{H}\bigr|^{\alpha}
\bigr\}.
\end{eqnarray*}
Since $\gamma>H-\frac{1}{2}$ and $\alpha>1-\frac{1}{2H}>H-\frac
{1}{2}$, we see that between jumps the paths $t\mapsto b(t,X_{t})$ are
H\"older continuous of order $H-\frac{1}{2}+\varepsilon$ for some
$\varepsilon>0$.
By the same argument as in Section~\ref{sec4}, it can be checked that $\hP\{
v\in L^2([0,T])\}=1$. Using the estimates
\begin{eqnarray*}
\bigl|b(t,X_{t})\bigr|\leq C\bigl(\bigl|b(0,x)\bigr|+t^{\gamma}+|X_{t}-x|^{\alpha}
\bigr)
\end{eqnarray*}
and
$\|X \|_{\infty}\leq C(1+|x|+\|B^{H}\|_{\infty}+|L|_{T})$, we deduce
that, for any $0\le r<t\le T$,
%
\begin{eqnarray}
\label{ualphart}\biggl | \int_{r}^{t}|u_{s}|\,\mathrm{d}s
\biggr|^{\alpha}&\leq& C\bigl(\bigl|b(0,x)\bigr|+t^{\gamma}+|x|^{\alpha}+\|X
\|^{\alpha}_{\infty}\bigr)^{\alpha
}(t-r)^{\alpha} .
\end{eqnarray}
In particular, we have
%
\begin{eqnarray}
\label{ualpha0t} \biggl| \int_{0}^{t}|u_{s}|\,\mathrm{d}s
\biggr|^{\alpha}&\leq& C\bigl(\bigl|b(0,x)\bigr|+t^{\gamma}+|x|^{\alpha}+\|X
\|^{\alpha}_{\infty}\bigr)^{\alpha
}t^{\alpha}
\nonumber
\\
&\leq&C \bigl(1+\bigl|b(0,x)\bigr|+|t|^{\gamma}+|x|^{\alpha}+\bigl(|x|+
\bigl\|B^{H}\bigr\|_{\infty
}+|L|_{T}\bigr)\bigr)^{\alpha}T^{\alpha}
\\
&\leq& C \bigl[1+\bigl|b(0,x)\bigr|^{\alpha}+t^{\alpha\gamma}+|x|^{\alpha}+\bigl\|
B^{H}\bigr\|^{\alpha}_{\infty}+|L|^{\alpha}_{T}
\bigr].\nonumber
\end{eqnarray}
Furthermore, one can also check that, by applying (\ref{ualphart}) and
(\ref{ualpha0t}), respectively,
%
\begin{eqnarray}
\label{Ak2-5} \bigl|A^k(t)\bigr| &\leq& C_1^H\max
_{1\leq i \leq k}\biggl|b\biggl(t,B_{t}^{H}+\int
_{0}^{t}u_s\,\mathrm{d}s-S_{i-1}+x
\biggr)-b\biggl(t,B_{t}^{H}+\int_{0}^{t}u_s\,\mathrm{d}s+x
\biggr)\biggr|
\nonumber
\\
&&{}\times \Biggl|\sum_{i=1}^{k-1} \biggl[
\frac{1}{(t-\sigma_{i})^{H-1/2}}-\frac{1} {
(t-\sigma_{i-1})^{H-1/2}} \biggr]+\frac{1}{(t-\sigma
_{k-1})^{H-1/2}}\Biggr |
\nonumber
\\
&&{}+C_1^H t^{1/2-H} \biggl|b\biggl(t,B_{t}^{H}+
\int_{0}^{t}u_s\,\mathrm{d}s+x\biggr)\biggr |
\nonumber
\\
&\leq&C \biggl\{(t-\sigma_{k-1})^{1/2-H}|L|_{T}^{\alpha}\\
&&\hphantom{C \bigl\{}{}+t^{1/2-H}
\biggl(\bigl|b(0,x)\bigr|+|t|^{\gamma}+\bigl\|B^{H}\bigr\|_{\infty}^{\alpha}+
\biggl|\int_{0}^{t}|u_s|\,\mathrm{d}s
\biggr|^{\alpha} \biggr) \biggr\}
\nonumber
\\
&\leq& C \bigl\{(t-\sigma_{k-1})^{1/2-H}|L|_{T}^{\alpha}\nonumber\\
&&\hphantom{C \bigl\{}{}+
t^{1/2-H}\bigl\| B^{H}\bigr\|^{\alpha}_{\infty}+
t^{1/2-H}\bigl[1+|x|+\bigl|b(0,x)\bigr|+|t|^{\gamma
}+|L|_{T}^{\alpha}
\bigr] \bigr\}
\nonumber
\\
&\leq& C \bigl\{(t-\sigma_{k-1})^{1/2-H}|L|_{T}^{\alpha}+
t^{1/2-H}\bigl\| B^{H}\bigr\|^{\alpha}_{\infty}+t^{1/2-H}
\bigl(1+|x|+\bigl|b(0,x)\bigr|+|t|^{\gamma}\bigr) \bigr\}\nonumber
\end{eqnarray}
and
%
\begin{eqnarray}
\label{Bk2-5} \bigl|B_{1}^{k}(t)\bigr|&\leq&\max
_{1\leq i \leq k}\biggl|b\biggl(t,B_{t}^{H}+\int
_{0}^{t}u_s\,\mathrm{d}s-S_{i-1}+x
\biggr)\biggr|t^{1/2-H}+K t^{\gamma+1/2-H}
\nonumber
\\
&\leq&Ct^{1/2-H} \biggl\{\bigl|b(0,x)\bigr|+|t|^{\gamma}+\bigl\|B^{H}
\bigr\|_{\infty
}^{\alpha}+|L_{T}|^{\alpha}+ \biggl|\int
_{0}^{t}u_s\,\mathrm{d}s \biggr|^{\alpha
}
\biggr\}
\\
&\leq& Ct^{1/2-H} \bigl\{1+|x|+\bigl|b(0,x)\bigr|+|L|_{T}^{\alpha}+
\bigl\|B^{H}\bigr\| ^{\alpha}_{\infty}+|t|^{\gamma} \bigr\},\nonumber
\end{eqnarray}
%
\begin{eqnarray}
\bigl|B_{2}^{k}(t)\bigr|&\leq& t^{H-1/2}\int
_{0}^{t}\frac{|\int_{r}^{t}u_s\,\mathrm{d}s+|B_{t}^{H}-B_{r}^{H}||^{\alpha
}}{(t-r)^{1/2+H}}r^{1/2-H}\,\mathrm{d}r
\nonumber
\\
&\leq&t^{H-1/2} \int_{0}^{t}
\frac{|\int_{r}^{t}u_s\,\mathrm{d}s|^{\alpha
}}{(t-r)^{1/2+H}}r^{1/2-H}\,\mathrm{d}r+t^{H-1/2} \int_{0}^{t}
\frac{|B_{t}^{H}-B_{r}^{H}|^{\alpha
}}{(t-r)^{1/2+H}}r^{1/2-H}\,\mathrm{d}r
\nonumber
\\
&\leq&t^{H-1/2}C\bigl(1+\bigl|b(0,x)\bigr|+|x|^{\alpha}+\|X
\|^{\alpha}_{\infty
}\bigr)^{\alpha}\nonumber\\[-8pt]\\[-8pt]
&&{}\times\int_{0}^{t}
\frac{(t-r)^{\alpha
}r^{1/2-H}}{(t-r)^{1/2+H}}\,\mathrm{d}r+C t^{1/2-H+\alpha(H-\varepsilon
)}G^{\alpha}\nonumber
\\
&\leq& t^{\alpha+H-1/2}C\bigl(1+\bigl|b(0,x)\bigr|+|x|^{\alpha}+\bigl\|B^{H}
\bigr\|_{\infty
}^{\alpha}+|L|_{T}^{\alpha}\bigr) +C
t^{1/2-H+\alpha(H-\varepsilon
)}G^{\alpha}
\nonumber
\\
&\leq& t^{\alpha+1/2-H}C \bigl\{1+|x|+\bigl|b(0,x)\bigr|+\bigl\|B^{H}
\bigr\|_{\infty}^\a +|L|^{\alpha}_{T} \bigr\}+C
t^{1/2-H+\alpha(H-\varepsilon)}G^{\alpha}.\nonumber
\end{eqnarray}
%
We can follow the same arguments of Lemma~\ref{Girsanov1} to show that
$v$ also satisfies the Novikov
condition (\ref{Novikov}), proving the theorem.
\end{pf}

Next, we show that the pathwise uniqueness holds for solutions to (\ref{SDE1}).
The proof is more or less standard, see \cite{RW} or \cite{S}, we
provide a sketch for completeness.

\begin{theorem}
\label{Uweak2}
Suppose that Assumption~\ref{Assump1} holds. Then two weak solutions
of SDE (\ref{SDE1}) defined on the same filtered probability space
with the same driving fBM $B^H$ and Poisson point process $L$
must coincide almost surely on their common existence interval.
\end{theorem}

\begin{pf} Let $X^{1}$ and $X^{2}$ be two weak solutions defined on
the same filtered probability space with the same driving $B^{H}$ and
$L$. Define $Y^+\stackrel{\triangle}{=}X^1\vee X^2$, and
$Y^-\stackrel{\triangle}{=}X^1\wedge X^2$. One
shows that both $Y^+$ and $Y^-$ both satisfy (\ref{SDE1}). In fact,
note that $X^1-X^2$ involves only Lebesgue integral, the occupation
density formula yields that the local time of $X^{1}-X^{2}$ at $0$ is
identically zero. Thus,
by Tanaka's formula,
\begin{eqnarray*}
\bigl(X^{1}_{t}-X^{2}_{t}
\bigr)^{+}=\int_{0}^{t}\bigl(b
\bigl(s,X_{s}^{1}\bigr)-b\bigl(s,X_{s}^{2}
\bigr)\bigr)I_{\{X^{1}_{s}-X^{2}_{s}>0\}}\,\mathrm{d}s.
\end{eqnarray*}
Then, note that $ Y^+=X^{2}+(X^{1}-X^{2})^{+}$, we have
\begin{eqnarray*}
Y^+_{t}&=&x+\int_{0}^{t}b
\bigl(s,X^{2}_{s}\bigr)\,\mathrm{d}s+B_{t}^{H}-L_{t}+
\int_{0}^{t}\bigl(b\bigl(s,X_{s}^{1}
\bigr)-b\bigl(s,X_{s}^{2}\bigr)\bigr)I_{\{X^{1}_{s}-X^{2}_{s}>0\}}\,\mathrm{d}s
\\
&=&x+\int_{0}^{t}b\bigl(s,X^{1}_{s}
\bigr)I_{\{X^{1}_{s}-X^{2}_{s}>0\}}\,\mathrm{d}s+\int_{0}^{t}b
\bigl(s,X^{2}_{s}\bigr)I_{\{X^{1}_{s}-X^{2}_{s}\leq0\}
}\,\mathrm{d}s+B_{t}^{H}-L_{t}
\\
&=&x+\int_{0}^{t}b\bigl(s, Y^+_{s}
\bigr)\,\mathrm{d}s+B_{t}^{H}-L_{t}.
\end{eqnarray*}
Similarly one shows that $Y^-_{t}$ satisfies SDE (\ref{SDE1}) as well.
We claim that
%
\begin{eqnarray}
\label{supX12=0} \hP \Bigl\{\sup_{0\leq t \leq T}\bigl(Y^+_{t}-Y^-_{t}
\bigr)=0 \Bigr\}=1.
\end{eqnarray}
Indeed, if $\hP \{\sup_{0\leq t \leq T}(Y^+_{t}-Y^-_{t})>0 \}>0$,
then there exists a rational number $r$ and $t>0$ such that
$\hP(Y^+_{t}>r>Y^-_{t})>0$. Since $\{Y^+_{t}>r\}=\{Y^-_{t}>r\}\cup\{
Y^+_{t}>r\geq Y^-_{t}\}$, we have
\[
\hP\bigl(Y^+_{t}>r\bigr)=\hP\bigl(Y^-_{t}>r\bigr)+\hP
\bigl(Y^+_{t}>r\geq Y^-_{t}\bigr)>\hP\bigl(Y^-_{t}>r
\bigr).
\]
This contradicts with the fact that $Y^+_{t}$ and $Y^-_{t}$ have the
same law, thanks to Theorem~\ref{Uweak1}. Thus,
(\ref{supX12=0}) holds, and consequently, $X^1\equiv X^2$, $\hP$-a.s.,
proving the theorem.
\end{pf}

\section{Existence of strong solutions}\label{sec7}

\setcounter{equation}{0}

Having proved the existence of the weak solution and pathwise
uniqueness, it is rather tempting to invoke the well-known
Yamada--Watanabe Theorem to conclude the existence of the strong
solution. However,
there seem to be some fundamental difficulties in the proof of such a
result, mainly because of
the lack of the independent increment property for an fBM, which is
crucial in the proof. It is also well known that, unlike an ODE, in the
case of stochastic differential equations, the existence of the strong
solution could be argued with assumptions on the coefficients being
much weaker than Lipschitz, due to the presence of the ``noise''. We
note that the argument in this section
is quite similar to \cite{GyongyPard} and \cite{NO}, with some
necessary adjustments for the presence of the jumps.

We begin by observing that the SDE (\ref{SDE1}) can be solved
pathwisely, as an ODE, when the coefficient $b$ is regular enough
(e.g., continuous in $(t,x)$, and uniformly Lipschitz in $x$).
Second, we claim that, under Assumption~\ref{Assump1} it suffices to
prove the existence of the strong solution when the coefficient $b$ is
uniformly bounded. Indeed, if we consider the following
family of SDEs:
%
\begin{eqnarray}
\label{SDER} X_{t}=x+\int_{0}^{t}{b}_{R}(s,X_{s})\,\mathrm{d}s+B_{t}^{H}-L_{t},\qquad
  t\in [0,T], R>0,
\end{eqnarray}
where $b_R$ is the truncated version of $b$:
$b_{R}(t,x)=b(t,(x\wedge R)\vee(-R))$, $(t,x)\in[0,T]\times\hR$,
then for each $R$, $b_R$ is bounded, hence (\ref{SDER}) has a strong
solution, denoted by $X^R$, defined on $[0,T]$, and we can now assume
that they all live on a common probability space. Now note that for
$R_1<R_2$, one has $b_{R_1}\equiv b_{R_2}$ whenever $|x|\le R_1$, thus
by the pathwise uniqueness, it is easy to see that $X_t^{R_1}\equiv
X_t^{R_2}$, for $t\in[0,\t_{R_1}]$, $\hP$-a.s., where
$\tau_{R}\stackrel{\triangle}{=}\inf\{t>0\dvtx  |X_{t}^{R}|\geq R\}
\wedge T$. Therefore, we
can almost surely extend the solution to $[0,\t)$, where $\t\stackrel
{\triangle}{=}
\lim_{R\to
\infty}\t_R$. Furthermore, it was shown (see, e.g., (\ref
{Xinfinity})) that $X$ will never explode
on $[0,\t)$. Consequently, we must have $\tau=T$, $\hP$-a.s.

We now give our main result of this section.

\begin{theorem}
\label{Ustrong}
Assume that $b(t,x)$ satisfies Assumption~\ref{Assump1}. Then there
exists a unique strong solution SDE (\ref{SDE1}).
\end{theorem}

%

%


%

The proof of Theorem~\ref{Ustrong} follows an argument by Gy\"
ongy and Pardoux \cite{GyongyPard}, using
the so-called Krylov estimate (cf. \cite{Krylov}). We note that by the
argument preceding the theorem
we need only consider the case when the coefficient $b$ is bounded. The
following lemma is thus crucial.

\begin{lem}
\label{krylov1}
Suppose that the coefficient $b$ satisfies Assumption~\ref{Assump1}
and is uniformly bounded by a
constant $C>0$. Suppose also that $X$ is a strong solution to SDE (\ref
{SDE1}). Then, there exist $\beta>1$ and $\zeta>1+H$ such that for
any measurable nonnegative function $g\dvtx  [0,T]\times\hR\mapsto\hR
_+$, it holds that
%
\begin{eqnarray}
\label{krylovI} \hE\int_{0}^{T}g(t,X_{t})\,\mathrm{d}t
\leq M \biggl(\int_{0}^{T}\int
_{\hR
}g^{\beta\zeta}(t,x)\,\mathrm{d}x\,\mathrm{d}t \biggr)^{1/\beta\zeta},
\end{eqnarray}
where $M$ is a constant defined by
%
\begin{eqnarray}
\label{G} M\stackrel{\triangle} {=}J^{1/\zeta^{\prime}\beta} F^{1/\a},
\end{eqnarray}
in which
%
\begin{eqnarray}
\label{KJ} F\stackrel{\triangle} {=} \biggl\{\tilde{\hE}\exp \biggl\{2
\a^2\int_0^T v_t^2\,\mathrm{d}t
\biggr\} \biggr\}^{1/2},\qquad    J\stackrel{\triangle} {=}
\frac{(2\pi)^{1/2-\zeta^{\prime}/2}
T^{1+(1-\zeta^{\prime})H}}{\sqrt {\zeta^{\prime}}(1+(1-\zeta^{\prime})H)}
\end{eqnarray}
and $\frac{1}{\alpha}+\frac{1}{\beta}=1$, $\frac{1}{\zeta}+\frac
{1}{\zeta^{\prime}}=1$.
\end{lem}

\begin{pf} Let $(\Omega, \cF, \hP;\hF)$ be a filtered probability
space on which are defined a fBM $B^H$, a Poisson point process $L$
of class (QL) and independent of $B^H$, and $X$ is the strong solution
to the corresponding SDE (\ref{SDE1}). Let $W$ be an $\hF$-Brownian
motion such that $B^H=\int_0^\cdot K_H(t,s)\,\mathrm{d}W_s$. Recall from (\ref
{v-2}) the process $v=K_{H}^{-1} (\int_{0}^{\cdot}b(r, X_r)\,\mathrm{d}r
)$, and define a new measure $\tilde{\hP} $ by
%
\begin{eqnarray}
\label{ZT} \frac{\mathrm{d}\tilde{\hP} }{\mathrm{d}\hP}\stackrel{\triangle} {=}\exp \biggl\{ -\int
_{0}^{T}v_{t} \,\mathrm{d}W_{t} -
\frac{1}{2}\int_{0}^{T}v_{t}^{2}\,\mathrm{d}t
\biggr\}\stackrel{\triangle } {=}Z_T^{-1}.
\end{eqnarray}
Then, in light of Lemmas \ref{Girsanov} and \ref{Girsanov1}, we know
that $\tilde\hP$ is a
probability measure
under which $\tilde{W}_{t}=W_{t}+ \int_{0}^{t}v_r\,\mathrm{d}r$ is a Brownian
motion, $\tilde{B}_{t}^{H}=\int_{0}^{t}K_{H}(t,s)\,\mathrm{d}\tilde{W}_{s}$ is
a fBM, and
$L$ remains a Poisson point process with same parameters and is
independent of $\tilde{B}^{H}$. Hence, under $\tilde{\hP}$,
$X_{t}=x+\tilde B^H_t-L_t$ has the density function:
%
\begin{eqnarray}
 p_{t}(y)=\int_{\hR}
\frac{1}{\sqrt{2\pi}t^{H}}\mathrm{e}^{-(y+z-x)^{2}/2
t^{2H}}f_L(t, z)\,\mathrm{d}z,
\end{eqnarray}
where $f_L(t,\cdot)$ is the density function of $L_t$.

Now, applying H\"{o}lder's inequality we have
%
\begin{eqnarray}
\label{EgX} \hE\int_{0}^{T}g(t,X_{t})\,\mathrm{d}t=
\tilde\hE \biggl\{Z_T\int_0^Tg(t,X_t)\,\mathrm{d}t
\biggr\}\le\bigl\{\tilde{\hE}\bigl[Z_T^\a\bigr]\bigr\}
^{1/\alpha} \biggl\{\tilde{\hE}\int
_{0}^{T}g^{\beta}(t,X_{t})\,\mathrm{d}t
\biggr\}^{1/\beta},
\end{eqnarray}
where $1/\alpha+1/\beta=1$. Rewriting $v_{t}$ as
$v_{t}=K_{H}^{-1} (\int_{0}^{\cdot}b(r, \tilde
{B}_{r}^{H}-L_{r}+x )\,\mathrm{d}r )(t)$,
we can follow the same argument as the proof of Lemmas \ref{Girsanov} and \ref{Girsanov1} to get,
$\tilde{\hE}\mathrm{e}^{2 \alpha^{2}\int_{0}^{T}v_{t}^{2}\,\mathrm{d}t}<\infty$.
Therefore, $\exp \{2\alpha\int_{0}^{t}v_{s} \,\mathrm{d}\tilde{W}_{s} -
2\alpha^{2}\int_{0}^{t}v^{2}_{s} \,\mathrm{d}s \}$ is a $\tilde{\hP
}$-martingale, and consequently, applying H\"older's inequality we obtain
%
\begin{eqnarray}
\label{EZT} \tilde\hE\bigl[Z^\a_T\bigr]&=&\tilde{
\hE}\exp \biggl\{\alpha\int_{0}^{T}v_{t}
\,\mathrm{d}W_{t} +\frac{\alpha}{2}\int_{0}^{T}v_{t}^{2}\,\mathrm{d}t
\biggr\}
\nonumber
\\
&=&\tilde{\hE}\exp \biggl\{\alpha\int_{0}^{T}v_{t}
\,\mathrm{d}\tilde{W}_{t} -\frac{\alpha}{2}\int_{0}^{T}v_{t}^{2}\,\mathrm{d}t
\biggr\}\nonumber
\\
&=& \tilde{\hE}\exp \biggl\{\alpha\int_{0}^{T}v_{t}
\,\mathrm{d}\tilde{W}_{t} -\alpha^{2}\int_{0}^{T}v^{2}_{t}
\,\mathrm{d}t+ \biggl(\alpha^{2}-\frac{\alpha
}{2}\biggr)\int
_{0}^{T}v_{t}^{2}\,\mathrm{d}t \biggr
\}
\\
&\leq& \biggl(\tilde{\hE}\exp \biggl\{2\alpha\int_{0}^{T}v_{t}
\,\mathrm{d}\tilde{W}_{t} - 2\alpha^{2}\int_{0}^{T}v^{2}_{t}
\,\mathrm{d}t \biggr\} \biggr)^{1/2} \biggl(\tilde{\hE}\exp \biggl\{\bigl(2
\alpha^{2} - \alpha\bigr)\int_{0}^{T}v_{t}^{2}\,\mathrm{d}t
\biggr\} \biggr)^{1/2}
\nonumber
\\
&\le& \biggl(\tilde{\hE}\exp \biggl\{2\alpha^{2} \int
_{0}^{T}v_{t}^{2}\,\mathrm{d}t \biggr
\} \biggr)^{1/2}<\infty.
\nonumber
\end{eqnarray}
%
On the other hand, applying H\"older's
inequality with $1/\zeta+1/\zeta^{\prime}=1$, $\zeta>H+1$ yields
%
\begin{eqnarray}
\label{tEgb} \tilde{\hE}\int_{0}^{T}g^{\beta}(t,X_{t})\,\mathrm{d}t&=&
\int_{0}^{T} \int_{\hR}g^{\beta}(t,y)p_{t}(y)\,\mathrm{d}y
\,\mathrm{d}t \nonumber\\[-8pt]\\[-8pt]
&\leq& \bigl\|g^\beta
\bigr\|_{L^{\zeta}([0,T]\times\hR)} \bigl\|p_\cdot(\cdot)\bigr\| _{L^{\zeta^{\prime}}([0,T]\times\hR)}.\nonumber 
\end{eqnarray}

Now, by the generalized Minkowski inequality (cf., e.g., \cite
{SKM}, (1.33)), we have
%
\begin{eqnarray}
\label{pty} \int_{\hR}\bigl[p_t(y)
\bigr]^{\gamma^{\prime}}\,\mathrm{d}y&=&\int_{\hR} \biggl\{\int
_{\hR
}\frac{1}{\sqrt{2\pi}t^{H}}\mathrm{e}^{-(y+z-x)^{2}/2 t^{2H}} 
f_L(t,z)\,\mathrm{d}z \biggr\}^{\zeta^{\prime}}\mathrm{d}y
\nonumber
\\
&\leq& \biggl\{\int_{\hR} \biggl[\int_{\hR}
\biggl(\frac{1}{\sqrt {2\pi}t^{H}}\mathrm{e}^{-(y+z-x)^{2}/2 t^{2H}} 
f_L(t,z)
\biggr)^{\zeta^{\prime}}\,\mathrm{d}y \biggr]^{1/\zeta^{\prime}}\,\mathrm{d}z \biggr\}^{\zeta
^{\prime}}
\\
&=& \biggl\{\int_{\hR} f_L(t,z) 
\biggl[\int_{\hR} \biggl(\frac{1}{\sqrt{2\pi}t^{H}}\mathrm{e}^{-(y+z-x)^{2}/2
t^{2H}}
\biggr)^{\zeta^{\prime}} \,\mathrm{d}y \biggr]^{1/\zeta^{\prime}}\,\mathrm{d}z \biggr\}^{\zeta
^{\prime}}.
\nonumber
\end{eqnarray}
The direct calculation gives
\[
\int_{\hR} \biggl(\frac{1}{\sqrt{2\pi}t^{H}}\mathrm{e}^{-(y+z-x)^{2}/2
t^{2H}}
\biggr)^{\zeta^{\prime}} \,\mathrm{d}y=(2\pi)^{1/2-\zeta^{\prime}/2}\bigl(\zeta ^{\prime}
\bigr)^{-1/2}t^{(1-\zeta^{\prime})H}.
\]
Plugging this into (\ref{pty}), we obtain
\begin{eqnarray*}
\int_{\hR}\bigl[p_t(y)\bigr]^{\zeta^{\prime}}\,\mathrm{d}y
&=&(2\pi)^{1/2-\zeta^{\prime}/2}\bigl(
\zeta^{\prime}\bigr)^{-1/2}t^{(1-\zeta^{\prime})H} \biggl(\int
_{\hR} 
f_L(t,z)\,\mathrm{d}z
\biggr)^{\zeta^{\prime}}
\\
&=&(2\pi)^{1/2-\zeta^{\prime}/2}\bigl(\zeta^{\prime}\bigr)^{-1/2}t^{(1-\sigma^{\prime})H}.
\end{eqnarray*}
Since $\zeta>H+1$, this leads to that
%
\begin{eqnarray}
\label{ptyg'} \bigl\|p_\cdot(\cdot)\bigr\|_{L^{\zeta^{\prime}}([0,T]\times\hR)}\le
J^{1/\zeta^{\prime}},
\end{eqnarray}
where $J$ is defined by (\ref{KJ}). Finally, noting that $\|g^\beta\|
^{1/\beta}_{L^{\zeta}([0,T]\times\hR)}
=\|g\|_{L^{\beta\zeta}([0,T]\times\hR)}$, the estimate (\ref
{krylovI}) then follows from (\ref{EgX}), (\ref{EZT}), (\ref{tEgb}),
and (\ref{ptyg'}).
\end{pf}
\begin{pf*}{Proof of Theorem~\ref{Ustrong}} Since the proof
is more or less standard, we only give a sketch for the completeness.
We refer to \cite{Krylov}, \cite{GyongyPard} and/or \cite{NO} for
more details.

We need only prove the existence. We assume that the coefficient $b$ is
bounded (by $C>0$) and satisfies Assumption~\ref{Assump1}. Let $\{
b_n(\cdot,\cdot)\}_{n=1}^\infty$ be a sequence of the mollifiers
of $b$, so that all $b_n$'s are smooth, have the
same bound $C$, and satisfy Assumption~\ref{Assump1} with the same parameters.

Next, for $n\leq k$ we define
$\tilde{b}_{n,k}\stackrel{\triangle}{=}\bigwedge_{j=n}^{k}b_{j}$
and $\tilde
{b}_{n}\stackrel{\triangle}{=}\bigwedge_{j=n}^{\infty}b_{j}$. Then
clearly, each
$\tilde{b}_{n,k}$ is continuous, and uniformly Lipschitz in $x$,
uniformly with respect to $t$. Furthermore, it holds that
\[
\tilde{b}_{n,k}\downarrow\tilde{b}_{n}, \qquad \mbox{as } k
\rightarrow \infty, \qquad   \tilde{b}_{n}\uparrow b, \qquad \mbox{as } n
\rightarrow\infty,
\]
for almost all $x$. Now for fixed $n$, $k$, consider SDE
%
\begin{eqnarray}
\label{SDEn} X_{t}=x+\int_{0}^{t}
\tilde{b}_{n,k}(s,X_{s})\,\mathrm{d}s+B_{t}^{H}-L_{t},\qquad
  t\ge0.
\end{eqnarray}
As a pathwise ODE, (\ref{SDEn}) has a unique strong solution $\tilde
{X}^{n,k}$, and comparison theorem holds, that is,
$\{\tilde{X}^{n,k}\}$
decrease with $k$. Furthermore, since $\tilde b_{n,k}$'s are uniformly
bounded by $C$, the solutions
$\tilde X^{n,k}$ are pathwisely uniformly bounded, uniformly in $n$ and
$k$. Thus $X^n_t\stackrel{\triangle}{=}\lim_{k\to
\infty} \tilde X^{n,k}_t$ exists, for all $t\in[0,T]$, $\hP$-a.s.
Since $b_n$'s are still Lipschitz, the standard stability result of ODE
then implies that
$\tilde{X}^{n}$ solves
\[
X_{t}=x+\int_{0}^{t}
\tilde{b}_{n}(s,X_{s})\,\mathrm{d}s+B_{t}^{H}-L_{t},\qquad
  t\in[0,T].
\]
Furthermore, the Dominated Convergence theorem leads to that the
estimate (\ref{krylovI}) holds for all $X^n$'s, for any bounded
measurable function $g$.

Next, since $\tilde{X}^{n,k}\leq\tilde{X}^{m,k}$, for $n\leq m \leq
k$, we see that $\tilde{X}_{n}$
increases as $n$ increases, thus $\tilde{X}^{n}$ converges, $\hP
$-almost surely, to some process $X$. The main task remaining is to
show that $X$ solves SDE (\ref{SDE1}), as $b$ is no longer Lipschitz.
In other words, we shall prove that
%
\begin{eqnarray}
\label{bntob} \lim_{n\to\infty}\hE\int_{0}^{T}\bigl|
\tilde b_{n}\bigl(t,X_{t}^n\bigr)-b(t,X_{t})\bigr|\,\mathrm{d}t=0.
\end{eqnarray}

To see this, we first note that
%
\begin{eqnarray}
\label{In} \hE\int_{0}^{T}\bigl|\tilde
b_{n}\bigl(t,X_{t}^n\bigr)-b(t,X_{t})\bigr|\,\mathrm{d}s
\leq I^n_{1}+I^n_2,
\end{eqnarray}
where
%
\begin{eqnarray}
\label{In12} I^n_{1}&\stackrel{\triangle} {=}&\sup
_{k}\hE\int_{0}^{T}\bigl|\tilde
b_{k}\bigl(t,X_{t}^n\bigr)-\tilde
b_{k}(t,X_{t})\bigr|\,\mathrm{d}t, \nonumber\\[-8pt]\\[-8pt]
  I^n_{2}
&\stackrel{\triangle} {=}&\hE\int_{0}^{T}\bigl|\tilde
b_{n}(t,X_{t})-b(t,X_{t})\bigr|\,\mathrm{d}t.\nonumber
\end{eqnarray}
Let $\kappa\dvtx \hR\rightarrow\hR$ be a smooth truncation function
satisfying $0\leq\kappa(z) \leq1$ for every $z$, $\kappa(z)=0$ for
$|z|\geq1$ and $\kappa(0)=1$. Then by Bounded Convergence theorem one has
%
\begin{eqnarray}
\label{kappato0} \lim_{R\to\infty}\hE\int_0^T
\bigl(1-\kappa(X_t/R)\bigr)\,\mathrm{d}t=0.
\end{eqnarray}
Now for any $R>0$, we apply Lemma~\ref{krylov1} with $\beta\zeta=2$
and note that both $\tilde b_n$ and $b$ are bounded by $C$ to get
\begin{eqnarray}
\label{I2n} I^n_2&=& \hE\int_{0}^{T}
\kappa(X_{t}/R)\bigl|\tilde b_{n}(t,X_{t}) -
b(t,X_{t})\bigr|\,\mathrm{d}t\nonumber\\
&&{}+ \hE\int_{0}^{T}
\bigl(1 - \kappa(X_{t}/R)\bigr)\bigl|\tilde b_{n}(t,X_{t})-b(t,X_{t})\bigr|\,\mathrm{d}t
\\
&\leq& M \biggl(\int_{0}^{T}\int
_{-R}^{R}\bigl|\tilde b_{n}(t,x)-b(t,x)\bigr|^{2}\,\mathrm{d}x\,\mathrm{d}t
\biggr)^{1/2}+2C\hE\int_{0}^{T}
\bigl(1-\kappa (X_{t}/R)\bigr)\,\mathrm{d}t.\nonumber
\end{eqnarray}
First letting $n\rightarrow\infty$ and then letting $R\rightarrow
\infty$, we get $\lim_{n\rightarrow\infty}I^n_{2}=0$.

To show that $\lim_{n\to\infty}I^n_1 =0$, we first note that by
(\ref{kappato0}), for any $\varepsilon>0$,
there exists $R_0>0$ such that
%
\begin{eqnarray}
\label{Kappalee} \hE\int_{0}^{T}\bigl|1-
\kappa(X_{t}/R_0)\bigr|\,\mathrm{d}t<\varepsilon.
\end{eqnarray}
Second, since $\{b_{n}\}$ converge to $b$ almost everywhere, the
Bounded Convergence theorem then shows that
$\tilde b_n$ converges to $b$ in $L^2_{T,R_0}\stackrel{\triangle
}{=}L^{2}([0,T]\times
[-R_0,R_0])$, hence $\{b_n, b\}_{n\ge1}$ is a compact set in
$L^2_{T,R_0}$. Thus, we can find finitely many bounded smooth function
$H_{1},\ldots, H_{N}$ such that
for each $k$, there is a $H_{i_k}$ so that
%
\begin{eqnarray}
\label{bkHi} \biggl(\int_{0}^{T}\int
_{-R_0}^{R_0}\bigl|\tilde b_{k}(t,x)-H_{i_k}(t,x)\bigr|^{2}\,\mathrm{d}r\,\mathrm{d}t
\biggr)^{1/2}<\e.
\end{eqnarray}
Now, we write
\begin{eqnarray*}
I^n_1=\sup_{k}\hE\int
_{0}^{T}\bigl|\tilde b_{k}
\bigl(t,X_{t}^n\bigr)-\tilde b_{k}(t,X_{t})\bigr|\,\mathrm{d}t
\leq\sup_{k} I_{1}(n,k)+I_{2}(n)+ \sup
_{k} I_{3}(k),
\end{eqnarray*}
where
\begin{eqnarray*}
\lleft\{ %
\begin{array} {l} I_{1}(n,k)=\hE\displaystyle \int
_{0}^{T}\bigl|\tilde b_{k}
\bigl(t,X_{t}^n\bigr)-H_{i_k}
\bigl(t,X_{t}^n\bigr)\bigr|\,\mathrm{d}t;
\\
I_{2}(n)=\displaystyle \sum_{j=1}^{N}\hE
\int_{0}^{T}\bigl|H_{j}
\bigl(t,X_{t}^n\bigr)-H_{j}(t,X_{t})\bigr|\,\mathrm{d}t;
\\
I_{3}(k)=\hE\displaystyle \int_{0}^{T}\bigl|\tilde
b_{k}(t,X_{t})-H_{i_k}(t,X_{t})\bigr|\,\mathrm{d}t.
\end{array} %
\rright.
\end{eqnarray*}

It is obvious that $\lim_{n\rightarrow\infty}I_{2}(n)=0$.
Furthermore, since the estimate (\ref{krylovI}) holds with $\beta
\zeta=2$ for
all $X^n$'s, similar to (\ref{I2n}) we have
\begin{eqnarray*}
I_{1}(n,k) 
&\leq& M
\biggl( \int_{0}^{T}\int_{-R_0}^{R_0}\bigl|
\tilde b_{k}(t,x)-H_{i_k}(t,x)\bigr|^{2}\,\mathrm{d}x\,\mathrm{d}t
\biggr)^{1/2}+C_{1} \hE\int_{0}^{T}
\bigl(1-\kappa\bigl(X_{t}^{(n)}/R_0\bigr)
\bigr)\,\mathrm{d}t,
\end{eqnarray*}
where $C_{1}$ is a constant depending on $C$ and
$\max_{1\leq i \leq N}\|H_{i}\|_{\infty}$.
Hence, by (\ref{Kappalee}), (\ref{bkHi}), and the Dominated
Convergence theorem again we have
\begin{eqnarray*}
\lim_{n\rightarrow\infty}\sup_{k}I_{1}(n,k)
\leq M \varepsilon +C_{1} \hE\int_{0}^{T}
\bigl(1-\kappa(X_{t}/R_0)\bigr)\,\mathrm{d}t\leq (M+C_{1})
\varepsilon.
\end{eqnarray*}
Similarly, we have $\sup_{k}I_{3}(k)\leq(M+C_{1})\varepsilon$.
Letting $\e\to0$ we obtain
$\lim_{n\rightarrow\infty}I^n_1=0$. The proof is now complete.
\end{pf*}

\section*{Acknowledgments}
The authors would like to thank referees for their helpful and valuable
suggestions. Lihua Bai is supported in part by NNSF of China grants \#11171164 and
\#11001136. Part of this work was completed while this author was
visiting Department of Mathematics,
University of Southern California, whose hospitality is greatly
appreciated. Jin Ma is supported in part by NSF grants
\#0806017 and \#1106853.




\printhistory

\end{document}